\documentclass[11pt,a4paper]{article}

\usepackage[margin=2.7cm]{geometry}
\usepackage[T1]{fontenc}
\usepackage[utf8]{inputenc}
\usepackage{lmodern}
\usepackage{microtype}
\usepackage{amsmath,amsthm,amssymb,amsfonts,mathtools}
\usepackage[numbers,sort&compress]{natbib}
\usepackage{graphicx}
\usepackage{enumitem}
\usepackage{aliascnt}
\usepackage[colorlinks,linkcolor=blue,citecolor=blue,urlcolor=blue]{hyperref}
\hypersetup{pdftitle={Characterizing M-estimators using canonical divergences from information geometry},pdfauthor={Paul Navas}}

\numberwithin{equation}{section}

\theoremstyle{plain}
\newtheorem{theorem}{Theorem}[section]
\newaliascnt{lemma}{theorem}
\newtheorem{lemma}[lemma]{Lemma}
\aliascntresetthe{lemma}
\newaliascnt{proposition}{theorem}
\newtheorem{proposition}[proposition]{Proposition}
\aliascntresetthe{proposition}
\newaliascnt{corollary}{theorem}
\newtheorem{corollary}[corollary]{Corollary}
\aliascntresetthe{corollary}
\theoremstyle{definition}
\newaliascnt{definition}{theorem}
\newtheorem{definition}[definition]{Definition}
\aliascntresetthe{definition}
\newaliascnt{assumption}{theorem}
\newtheorem{assumption}[assumption]{Assumption}
\aliascntresetthe{assumption}
\newaliascnt{remark}{theorem}
\newtheorem{remark}[remark]{Remark}
\aliascntresetthe{remark}
\newtheorem*{example}{Example}

\DeclareMathOperator{\grad}{grad}
\DeclareMathOperator{\Hess}{Hess}

\newcommand{\E}{\mathbb{E}}
\newcommand{\R}{\mathbb{R}}
\newcommand{\simplex}{\Delta^{n-1}}
\newcommand{\csimplex}{\bar\Delta^{n-1}}

\newcommand{\T}{^{\top}}
\newcommand{\loss}{\ell}

\newcommand{\bary}{b}
\newcommand{\musflat}[1]{{#1}^{\flat}}
\newcommand{\mussharp}[1]{{#1}^{\sharp}}
\newcommand{\Ff}{\mathcal{F}}
\newcommand{\Mm}{\mathcal{M}}

\defcitealias{dimitriadis2024}{DFZ}

\AtBeginDocument{%
}

\begin{document}

\title{Characterizing M-estimators using canonical divergences from information geometry}
\author{Paul Navas\\[3pt]
  \normalsize Department of Statistics, TU Dortmund University\\
  \normalsize \texttt{paul.navas@tu-dortmund.de}}
\date{\today}
\maketitle

\begin{abstract}
On Euclidean space, the losses that produce consistent M-estimators of a functional are its strictly consistent scoring functions \citep{dimitriadis2024}, and Osband's principle characterizes these through identification functions, which define method-of-moments estimators. We extend this characterization to the probability simplex. We define the residual, which on Euclidean space is a subtraction, as the inverse exponential map of an affine connection. It identifies the exponential barycenter of the connection: the arithmetic, normalized geometric and Fr\'echet means for the mixture, exponential and Levi-Civita connections. Its integral along geodesics is the canonical divergence of \citet{ayamari2015}. A geometric form of Osband's principle decides which losses each connection admits. For the flat connections every strictly convex generator yields a consistent Bregman loss. For the curved $\alpha$-connections, half the squared Fisher-Rao distance is (up to scale) the only consistent loss for the Fr\'echet mean, and for $\alpha \notin \{-1, 0, 1\}$ no such loss exists and only the method of moments remains. The residual is integrable only at $\alpha = 0$, sharpening the integrability claim of Ay and Amari. The characterization yields the estimating equations of regression on the simplex and the calibration decomposition of compositional forecasts.
\end{abstract}



\section{Introduction}\label{sec:Intro}

An M-estimator returns the minimizer of an empirical risk over a model $M$,
\begin{equation}\label{eq:m-estimator}
\hat q = \arg\min_{q \in M} \frac{1}{N} \sum_{k=1}^{N} \loss(q, x_k).
\end{equation}
With squared loss it is consistent for the mean, with absolute loss for the median, and with squared loss on the logarithms for the geometric mean. Conversely, the losses that estimate the mean consistently are exactly the Bregman losses of strictly convex generators \citep{savage1971,gneiting2011}. \citet*{dimitriadis2024}, DFZ in what follows, proved that, under richness conditions, the losses that produce consistent M-estimators of a functional are exactly the strictly consistent scoring functions for that functional in the sense of forecast evaluation \citep{gneiting2011}. Characterizing M-estimators and characterizing consistent scoring functions are therefore one problem. Through Osband's principle, the identification functions of a functional characterize its strictly consistent scoring functions \citep{osband1985,fisslerziegel2016}. An identification function is a moment function, and setting its sample average to zero estimates the functional by the method of moments, with no loss involved. By Osband's principle the first-order condition of every strictly consistent loss is such a moment equation, with the identification function reweighted by a matrix that depends on the argument alone; characterizing the losses therefore means characterizing the moment equations that come from a minimization.

For a Bregman loss $D_\psi[x:q] = \psi(x) - \psi(q) - \nabla\psi(q) \cdot (x - q)$ the derivative with respect to the argument $q$ is
\[
\nabla_q D_\psi[x:q] = -\nabla^2\psi(q)\,(x - q),
\]
a matrix that depends on $q$ alone, applied to the residual $x - q$. By Osband's principle, the score of every consistent loss for the mean has this form. An M-estimator consists of five pieces: the sample space, the loss, the risk that averages it, the score that differentiates it, and the estimating equation that sets the average score to zero. On the probability simplex only the score changes in substance. The risk is still an average, linear in the distribution, the estimating equation still sets an average to zero, and the characterization determines the loss. The residual is a subtraction, and on $\R^k$ subtraction and the notion of straight line are canonical. On the simplex both depend on the chosen geometry, with one residual, and hence one score, for each geometry.

Information geometry \citep{amari2000,ayjost2017} describes a notion of straight line on a manifold by an affine connection, whose straight curves are its geodesics. The mixture, exponential and Levi-Civita connections of the simplex are presented in \citet[Section~2]{navasconvex}. The simplex is an open subset of the hyperplane $\sum_i p_i = 1$ of $\R^n$. In the geometry induced from $\R^n$ the difference $x - q$ is still defined, sums to zero, and is a direction at $q$. This induced geometry is that of the mixture connection, whose geodesics are the segments $(1 - t)q + tx$. The geodesics of the exponential connection are the curves that are straight in the logarithms of the probabilities, up to normalization. Those of the Levi-Civita connection of the Fisher-Rao metric are the curves that become great circles under $p \mapsto 2\sqrt{p}$. For each connection, the residual at $q$ is the initial velocity of the geodesic that starts at $q$ and reaches $x$ at time one, the inverse exponential map $X(q, x)$. On $\R^k$ it is $x - q$. The point at which the expected residual vanishes is the exponential barycenter of the distribution for that connection \citep{pennec2013}. The residual is thus a moment function for the barycenter, and the root of $\sum_k X(\hat q, x_k) = 0$ is its method-of-moments estimator. It is the arithmetic mean for the mixture connection and the normalized geometric mean for the exponential connection. For Fisher-Rao it is the Fr\'echet mean, whose estimator, the Levi-Civita-Fr\'echet estimator, is studied in \citet{navasconvex}.

On $\R^k$ the identity between score and residual, $\nabla_q \tfrac{1}{2}\|x - q\|^2 = -(x - q)$, comes from differentiating the loss. On the simplex the geometry determines the residual $X(q, x)$, and the identity becomes an equation whose solution determines the loss,
\[
\grad_q \loss(\cdot, x) = -X(q, x),
\]
where the gradient turns the derivative of $\loss$, a linear form on the directions at $q$, into a direction at $q$ through the chosen metric, for instance Fisher-Rao. Since the simplex is contractible, a solution exists exactly when the residual field, turned into a linear form by the metric, has symmetric partial derivatives. The solution is then the integral of the residual along any path from $q$ to $x$. Integrating along the geodesics of the connection is the construction by which \citet{ayamari2015} define the canonical divergence. With the Fisher-Rao metric it is the Kullback-Leibler divergence in its two orders for the mixture and exponential connections, and half the squared Fisher-Rao distance, $d^2/2$, for the Levi-Civita connection.

Osband's principle also has a geometric form. Every strictly consistent loss for the exponential barycenter of a connection satisfies
\[
d_q \bar\loss(\cdot, F) = -H_q(\bar X(q, F), \cdot)
\]
for a field $H$ of bilinear forms on the tangent spaces, the Osband tensor of the loss; $\bar\loss$ and $\bar X$ are the expected loss and the expected residual under $F$. When the canonical divergence is a consistent loss, its Osband tensor is the metric itself. For the flat connections, the mixture and the exponential one, the Osband tensors are the Hessians of convex generators in the affine coordinates of the connection, and the consistent losses are the Bregman divergences, as on $\R^k$. The main contribution of this paper is a rigidity result for the curved connections. The $\alpha$-connections of the Fisher-Rao metric \citep{amari2000} run from the mixture connection at $\alpha = -1$ through the Levi-Civita connection at $\alpha = 0$ to the exponential connection at $\alpha = 1$. For them, curvature forces the Osband tensor of every loss in the regularity class of Osband's principle to be a constant multiple of the metric whenever $\alpha^2 \neq 1$ and $n \geq 3$. At $\alpha = 0$ the only consistent loss in this class, up to scale and an additive function of the observation, is $d^2/2$, strictly consistent for the Fr\'echet mean of every distribution on the open simplex. For $\alpha \notin \{-1, 0, 1\}$ no consistent loss exists in this class. The barycenter equation $\sum_k X(\hat q, x_k) = 0$ then does not come from the first-order conditions of a minimization, but it still defines a method-of-moments estimator, whose asymptotic normality requires only identification (Section~\ref{sec:regression}).

The $\alpha$-divergences of the cone of positive measures, restricted to the simplex, elicit normalized power means in closed form. The cone is dually flat for every $\alpha$, and the estimator averages powers of the observations and normalizes. At $\alpha = 0$ this power mean and the Fr\'echet mean are the extrinsic and the intrinsic means of \citet{bhattacharya2003}. Both are elicitable, with twice the squared Hellinger distance and $d^2/2$ as losses, the choice between them also determining the chosen functional. The M-estimator of the $\alpha$-divergence of the cone and the method-of-moments estimator of the residual of the simplex's own $\alpha$-connection coincide at $\alpha = \pm 1$ and differ at $\alpha = 0$. The equivalence of \citetalias{dimitriadis2024} between consistency for a conditional functional and model consistency holds for regression models $m: \R^r \times \Theta \to \simplex$ with $\Theta \subseteq \R^s$, which include geodesic regression \citep{fletcher2013}. In such a model the estimating equation pairs the residual with the columns of the model differential
\[
d_\theta m(z, \theta): \R^s \to T_{m(z, \theta)}\simplex, \quad v \mapsto \sum_{j=1}^{s} v_j \frac{\partial m(z, \theta)}{\partial \theta_j},
\]
the directions in which the fitted composition moves when a parameter changes. These columns act as instruments, as the covariates do in the normal equations of linear regression. The estimating equation is thus a method-of-moments equation with instruments, as in the generalized method of moments of \citet{hansen1982}. It yields the asymptotic normality of the estimators, and for the Fr\'echet mean this holds for every distribution on the open simplex, without moment conditions. Any invertible transformation of the residual, a gauge transformation, leaves the barycenter unchanged but changes the estimating equation of the model. All the resulting equations are solved at the true parameter, but in general only the gauges of consistent losses turn them into first-order conditions of a minimization, and at $\alpha = 0$ these are the multiples of the identity gauge, that of $d^2/2$. The influence function of an M-estimator depends on the functional and not on the loss, and for the mean it is unbounded. In regression, where contaminated distributions lie outside the model, it depends on the loss, and on the simplex it is bounded in the response for every loss whose residual is bounded, the mixture, Levi-Civita and power-mean losses, as it is not for linear regression. For the Fr\'echet mean, whose residual has length below the diameter $\pi$, the gross-error sensitivity is finite for each distribution, although the bound we obtain is not uniform over distributions.

For compositional forecasts, the expected score of every consistent loss considered here decomposes into uncertainty, discrimination and miscalibration \citep{dimitriadisgneitingjordan2021}. At the dually flat ends $\alpha = \pm 1$ the miscalibration term is a divergence between the issued and the recalibrated forecast. For the power means it is this divergence multiplied by a factor that measures the concentration of the outcome, and for the Fr\'echet loss it lies between two multiples of the expected squared distance between the two forecasts, with constants fixed by the curvature and the diameter of the region containing forecasts and outcomes. The geometry also decides how forecasts are ranked. At the flat ends two consistent losses can rank the same pair of forecasters in opposite ways, as on the real line \citep{patton2020}. For the Fr\'echet mean with at least three categories every consistent loss yields the same ranking, so once the geometry is chosen, consistency alone decides which forecaster is better (Section~\ref{sec:ranking}).

\section{Elicitability, consistent losses and identification functions in Euclidean space}\label{sec:elicitability}

This section recalls, for completeness, the results of the works that serve as our main reference, \citet{gneiting2011}, \citet{fisslerziegel2016} and \citet*{dimitriadis2024,dfzidentification}.\footnote{The sign convention is adapted so that the residual is $x - q$, as in Section~\ref{sec:Intro}.}

Let $\Ff$ be a class of distributions on $\R^k$ and $T: \Ff \to \R^k$ a functional, such as the mean or a quantile. A loss $\loss(q, x)$ is evaluated at an argument $q \in \R^k$ and an observation $x$, and its expected value under $F$ is $\bar\loss(q, F) = \E_F\,\loss(q, X)$.

\begin{definition}\label{def:consistency}
A loss $\loss$ is strictly $\Ff$-consistent for $T$ if $\bar\loss(T(F), F) < \bar\loss(q, F)$ for every $F \in \Ff$ and every $q \neq T(F)$. A functional with a strictly consistent loss is elicitable on $\Ff$.
\end{definition}

Strict consistency makes the M-estimator \eqref{eq:m-estimator} consistent. By the law of large numbers the empirical risk converges to $\bar\loss(\cdot, F)$, whose unique minimizer is $T(F)$, and under the usual regularity conditions the minimizers converge as well.

Squared loss is strictly consistent for the mean on distributions with finite second moment, since $\bar\loss(q, F) = \E_F\|X - \E_F X\|^2 + \|q - \E_F X\|^2$. Absolute loss is consistent for the median, and the pinball loss for quantiles \citep{gneiting2011}. Not every functional is elicitable. Suppose $\loss$ is strictly consistent and $T(F_0) = T(F_1) = t$. For every mixture $F_\lambda = (1 - \lambda)F_0 + \lambda F_1$ in $\Ff$, the expected loss $\bar\loss(\cdot, F_\lambda) = (1 - \lambda)\bar\loss(\cdot, F_0) + \lambda\bar\loss(\cdot, F_1)$ is minimized at $t$, so $T(F_\lambda) = t$. The level sets of an elicitable functional are therefore convex \citep{osband1985}, a test the variance fails. The argument uses only that the expected loss is linear in the distribution, never linearity in the argument, and it holds unchanged on the simplex (Section~\ref{sec:barycenter}).

A map $V: \R^k \times \R^k \to \R^k$ is a strict identification function for $T$ on $\Ff$ if $\bar V(q, F) = \E_F V(q, X)$ vanishes exactly at $q = T(F)$. For the mean it is the residual, $V(q, x) = x - q$, and for the $\tau$-quantile it is $V(q, x) = \tau - \mathbf{1}\{x \le q\}$. In econometric terms $V$ is a moment function and $\bar V(q, F) = 0$ a moment condition, and the root of $\frac{1}{N}\sum_k V(q, x_k) = 0$ is the method-of-moments estimator of $T$, defined without any loss. Identification functions are unique up to a gauge: every strict identification function has the form $A(q)V(q, x)$ with $A(q)$ an invertible matrix \citep{dfzidentification}.\footnote{We call the passage from $V$ to $A(q)V$ a gauge transformation. The matrix $A(q)$ may depend on the argument but not on the observation, so it changes the identification function without changing where its expectation vanishes, and hence without changing the identified functional.} Osband's principle \citep{osband1985}, in the form of \citet[Theorem~3.2]{fisslerziegel2016}, connects identification and consistency. Suppose $T$ has a strict identification function $V$, and the class $\Ff$ and the expected loss satisfy their Assumptions (V1) and (S1). Then every strictly consistent loss satisfies
\begin{equation}\label{eq:osband-rk}
\nabla_q \bar\loss(q, F) = -h(q)\,\bar V(q, F)
\end{equation}
for a matrix $h(q)$ that depends on the argument alone. The first-order condition of every strictly consistent loss is thus the method-of-moments equation of the reweighted identification function $h(q)V$. For the loss $\tfrac{1}{2}\|x - q\|^2$ the matrix $h$ is the identity, and \eqref{eq:osband-rk} says that the gradient of the expected loss is the negative expected residual, $\nabla_q \bar\loss(q, F) = -(\E_F X - q)$. For other losses $h(q)$ reweights the expected residual before it becomes a gradient, and not every field of matrices can do so, since a gradient has a symmetric derivative.

For the mean, this symmetry, the second-order condition of \citet[Corollary~3.3]{fisslerziegel2016}, forces $h = \nabla^2\psi$ for a strictly convex $\psi$, and integrating \eqref{eq:osband-rk} along segments recovers the Bregman losses
\begin{equation}\label{eq:bregman-def}
\loss(q, x) = D_\psi[x:q] + a(x),
\end{equation}
with $a$ an integrable function of the observation alone \citep{savage1971,gneiting2011}. Their score, $-\nabla^2\psi(q)(x - q)$, is \eqref{eq:osband-rk} with $h = \nabla^2\psi$. Their expected loss splits into a dispersion term and a divergence from the mean \citep{banerjee2005},
\begin{equation}\label{eq:two-stage}
\E_F D_\psi[X:q] = \E_F D_\psi[X:\mu] + D_\psi[\mu:q], \quad \mu = \E_F X.
\end{equation}
Section~\ref{sec:calibration} uses this identity for the decomposition of expected scores into calibration and discrimination.

\section{Residuals on the simplex as inverse exponential maps}\label{sec:residuals}

The open simplex $\simplex = \{p \in \R^n : p_i > 0,\ \sum_i p_i = 1\}$ is the set of discrete distributions with full support on $n$ categories, or equivalently the set of compositions of $n$ positive parts. It is the slice $\sum_i m_i = 1$ of the open cone $\Mm = \{m \in \R^n : m_i > 0\}$ of positive measures, the space in which the ambient geometries of Section~\ref{sec:alpha} are defined. A tangent vector at $p$ is the velocity $v = \dot\gamma(0)$ of a curve $\gamma$ in the simplex through $p$. Differentiating $\sum_i \gamma_i(t) = 1$ yields
\begin{equation}\label{eq:tangent}
T_p\simplex = \{v \in \R^n : \textstyle\sum_i v_i = 0\},
\end{equation}
which is the same hyperplane at every $p$. Its elements are reallocations of probability mass. For $n = 3$, for instance, $v = (0.2, -0.15, -0.05)$ adds mass to the first category and takes it from the other two. The cone is open in $\R^n$, so $T_m\Mm = \R^n$.

A Riemannian metric assigns an inner product to each tangent space. The Fisher-Rao metric of the categorical family is
\begin{equation}\label{eq:fr-metric}
g_p(v, w) = \sum_{i=1}^{n} \frac{v_i w_i}{p_i}, \quad v, w \in T_p\simplex,
\end{equation}
and the same formula defines a metric on the cone, known there as the Shahshahani metric.\footnote{Other metrics can be used in place of $g$: a divergence induces a metric and a pair of dual connections through its second and third derivatives on the diagonal (\citealp{eguchi1992}; \citealp[Theorem~2.6]{eguchikomori2022}), which is how \citet[Section~2.3]{navasconvex} obtains the metrics of convex and exponentially concave potentials.} The weights $1/p_i$ make lengths depend on the base point: moving mass $0.01$ is a large change for a category of probability $0.02$ and a small one for a category of probability $0.5$. With $\xi_i = 2\sqrt{p_i}$ one has $d\xi_i = dp_i/\sqrt{p_i}$, so the round metric of the sphere, $\sum_i d\xi_i^2$, pulls back to $g$, and $\Pi: \simplex \to S_2^{n-1} \cap \R^n_{>0}$, $p \mapsto 2(\sqrt{p_1}, \dots, \sqrt{p_n})$, is an isometry onto the open positive orthant of the sphere of radius $2$ in $\R^n$, which has constant curvature $1/4$. The geodesic distance $d(p, q)$, the infimum of the lengths of curves in the simplex from $p$ to $q$, is therefore the length of the great-circle arc between $\Pi(p)$ and $\Pi(q)$,
\begin{equation}\label{eq:rao-distance}
d(p, q) = 2\arccos \sum_{i=1}^{n} \sqrt{p_i q_i}.
\end{equation}
Its supremum over the simplex, the diameter, equals $\pi$ and is not attained: $\sum_i \sqrt{p_i q_i}$ is positive on the open simplex and tends to $0$ when $p$ and $q$ tend to distinct vertices.

A connection $\nabla$ specifies how a vector field is differentiated along a curve. Its geodesics are the curves $\gamma$ with $\nabla_{\dot\gamma}\dot\gamma = 0$; in coordinates, $\ddot\gamma^k + \Gamma^k_{ij}(\gamma)\,\dot\gamma^i\dot\gamma^j = 0$, where the Christoffel symbols $\Gamma^k_{ij}$ represent the connection in the chart. A connection is flat when its Christoffel symbols vanish in some chart, and its geodesics are then the curves with coordinates $a + tb$ in that chart, for constant vectors $a$ and $b$. The cone and the simplex each have a mixture chart and an exponential chart. On the cone these are the identity and $\log: \Mm \to \R^n$, $m \mapsto (\log m_1, \dots, \log m_n)$. On the simplex they are $p \mapsto (p_1, \dots, p_{n-1})$ and the additive log-ratio $\operatorname{alr}: \simplex \to \R^{n-1}$, $p \mapsto (\log(p_1/p_n), \dots, \log(p_{n-1}/p_n))$, the natural parameter of the categorical exponential family \citep[Section~2.2]{navasconvex}. The mixture and exponential connections are the connections whose Christoffel symbols vanish in the respective chart, written $\nabla^{(m)}$, $\nabla^{(e)}$ on the cone and $\bar\nabla^{(m)}$, $\bar\nabla^{(e)}$ on the simplex. From $q \in \simplex$ with initial velocity $v \in T_q\simplex$ their geodesics in the simplex are
\[
\gamma^{(m)}(t) = q + tv, \quad \gamma^{(e)}(t)_i = \frac{q_i e^{t v_i/q_i}}{\sum_j q_j e^{t v_j/q_j}}.
\]
On the cone the mixture geodesics have the same form, and the exponential ones lack the denominator.

The Levi-Civita connection of $g$ is the unique connection with symmetric Christoffel symbols, $\Gamma^k_{ij} = \Gamma^k_{ji}$, that satisfies
\[
u\,g(v, w) = g(\nabla_u v, w) + g(v, \nabla_u w)
\]
for all vector fields $u, v, w$, where $u\,g(v, w)$ is the derivative of the function $g(v, w)$ along $u$. The geodesics for the Levi-Civita connection of the simplex are the (locally) shortest curves for $d$, the curves that $\Pi$ maps onto arcs of great circles. The mixture and exponential connections satisfy this identity with one connection in each term, $u\,g(v, w) = g(\nabla^{(m)}_u v, w) + g(v, \nabla^{(e)}_u w)$, on the cone and, with bars, on the simplex; they are dual with respect to $g$ \citep{amari2000}. The $\alpha$-connections
\begin{equation}\label{eq:alpha-connection}
\nabla^{(\alpha)} = \frac{1 - \alpha}{2}\nabla^{(m)} + \frac{1 + \alpha}{2}\nabla^{(e)}, \quad \alpha \in [-1, 1],
\end{equation}
interpolate between them, and the same formula with bars defines $\bar\nabla^{(\alpha)}$ on the simplex. By linearity $\nabla^{(\alpha)}$ and $\nabla^{(-\alpha)}$ are dual, and $\nabla^{(0)}$, dual to itself and with symmetric Christoffel symbols, is the Levi-Civita connection. \citet[Section~3.3]{ayamari2015} define $\bar\nabla^{(\alpha)}$ instead as the connection induced from the cone, $\bar\nabla^{(\alpha)}_u v = P(\nabla^{(\alpha)}_u v)$, where $P_q: \R^n \to T_q\simplex$, $w \mapsto w - (\sum_j w_j)\,q$, is the orthogonal projection with respect to $g_q$. The two definitions agree, since $P$ is linear and maps $\nabla^{(m)}$ and $\nabla^{(e)}$ to $\bar\nabla^{(m)}$ and $\bar\nabla^{(e)}$ \citep[Section~3.2]{ayamari2015}.

The cone is dually flat for every $\alpha$: $\nabla^{(\alpha)}$ and its dual $\nabla^{(-\alpha)}$ are both flat, since the $\alpha$-geodesics of the cone are affine in the chart $m \mapsto (m_1^{(1-\alpha)/2}, \dots, m_n^{(1-\alpha)/2})$ for $\alpha < 1$ \citep[Equation~(31)]{ayamari2015} and in the exponential chart for $\alpha = 1$. On the simplex only the ends of the family at $\alpha = \pm 1$ are flat. The difference of two connections is a tensor, so $K = \tfrac{1}{2}(\bar\nabla^{(e)} - \bar\nabla^{(m)})$ is one, and $\bar\nabla^{(\alpha)} = \bar\nabla^{(0)} + \alpha K$. The curvature tensor $\bar R^{(\alpha)}(u, v) = \bar\nabla^{(\alpha)}_u\bar\nabla^{(\alpha)}_v - \bar\nabla^{(\alpha)}_v\bar\nabla^{(\alpha)}_u - \bar\nabla^{(\alpha)}_{[u, v]}$ contains products of two covariant derivatives, so it is a polynomial of degree two in $\alpha$, $\bar R^{(\alpha)} = \bar R^{(0)} + \alpha A + \alpha^2 B$, with tensors $A$ and $B$ that do not depend on $\alpha$. Flatness at $\alpha = 1$ and $\alpha = -1$ implies $A = 0$ and $B = -\bar R^{(0)}$, and the constant curvature $1/4$ of the Fisher-Rao metric yields
\begin{equation}\label{eq:alpha-curvature}
\bar R^{(\alpha)}(u, v)w = \frac{1 - \alpha^2}{4}\,(g(v, w)\,u - g(u, w)\,v).
\end{equation}
For $n \ge 3$ and $\alpha^2 \ne 1$ the right-hand side does not vanish, and $\bar\nabla^{(\alpha)}$ has no flat chart.

The residual $x - q$ is the identification function of the mean, and by Osband's principle \eqref{eq:osband-rk} the expected score of every strictly consistent loss for the mean is $-h(q)$ times the expected residual. On $\R^k$ the residual has two descriptions. It is the difference between the observation and the argument, and it is the velocity of the straight line $t \mapsto q + t(x - q)$, which leaves $q$ at time zero and reaches $x$ at time one. On the simplex a difference of coordinates depends on the chart, but the second description only requires straight lines, and a connection provides them as its geodesics. We take the second description as the definition of the residual on the simplex, with one residual for each connection. The exponential map at $q$ sends a velocity $v$ to the point $\exp_q(v) = \gamma_{q,v}(1)$ reached at time one by the geodesic with $\gamma_{q,v}(0) = q$ and $\dot\gamma_{q,v}(0) = v$, for the velocities whose geodesic stays in the space up to that time. The residual is its inverse, $X(q, x) = \exp_q^{-1}(x)$, the velocity at $q$ of the geodesic that reaches $x$ at time one, defined where $\exp_q$ is invertible.

For a flat connection the first description returns in the flat chart $\varphi$. The geodesic from $q$ to $x$ is $t \mapsto \varphi^{-1}((1 - t)\varphi(q) + t\varphi(x))$, and differentiating at $t = 0$ one obtains
\begin{equation}\label{eq:flat-residual}
X(q, x) = (d\varphi_q)^{-1}(\varphi(x) - \varphi(q)),
\end{equation}
where $d\varphi_q$ is the differential of the chart at $q$, a linear isomorphism from the tangent space onto the space of coordinates. The residual is a subtraction in the flat chart, turned into a tangent vector at $q$ by the inverse differential of the chart. The images of the charts above are convex, so \eqref{eq:flat-residual} is defined for all $q$ and $x$. On the cone it reproduces the inverse exponential maps computed by \citet[Section~3]{ayamari2015},
\begin{equation}\label{eq:cone-residuals}
X^{(m)}(q, x) = x - q, \quad X^{(e)}(q, x)_i = q_i\log\frac{x_i}{q_i}, \quad X^{(\alpha)}(q, x)_i = \frac{2}{1 - \alpha}\,q_i\left(\left(\frac{x_i}{q_i}\right)^{(1-\alpha)/2} - 1\right),
\end{equation}
with $\alpha < 1$ in the last one, and on the simplex it produces $\bar X^{(m)}(q, x) = x - q$ and $\bar X^{(e)}(q, x)_i = q_i(\log(x_i/q_i) - \sum_j q_j\log(x_j/q_j))$. In both cases $\bar X = P_q X$: at the ends of the family, the residual of the simplex is the projected residual of the cone. For $\alpha^2 \neq 1$, i.e. where $\bar\nabla^{(\alpha)}$ is curved, there is no flat chart and the residual is obtained from the geodesics. At $\alpha = 0$ it is the Riemannian logarithm $\log_q(x)$ of \citet[Section~3.4]{navasconvex}, of length $d(q, x)$, and Section~\ref{sec:alpha} computes it for the other $\alpha$.

The differential of a function $f$ on the simplex at $q$ is the linear form, or covector, $d_qf: T_q\simplex \to \R$, $v \mapsto \frac{d}{dt}f(\gamma(t))|_{t=0}$, for a curve $\gamma$ in the simplex with $\gamma(0) = q$ and $\dot\gamma(0) = v$. The metric identifies linear forms with tangent vectors. Lowering, $w \mapsto \musflat{w} = g_q(w, \cdot)$, is a linear isomorphism from $T_q\simplex$ onto its dual, and its inverse, raising, is written $\omega \mapsto \mussharp{\omega}$. The gradient is the raised differential, $\grad_q f = \mussharp{(d_qf)}$, the tangent vector with $g_q(\grad_q f, v) = d_qf(v)$ for every $v \in T_q\simplex$. On the cone the Fisher-Rao gradient is $(\grad_m f)_i = m_i\,\partial_i f(m)$, and on the simplex the gradient is the projection $P_q$ of this vector \citep[Equations~(16) and (19)]{ayamari2015},
\begin{equation}\label{eq:grad-simplex}
(\grad_q f)_i = q_i\left(\partial_i f(q) - \sum_j q_j\,\partial_j f(q)\right),
\end{equation}
for any extension of $f$ to a neighbourhood of the simplex in $\R^n$. The vector of partial derivatives is not a tangent vector of the simplex, since its components need not sum to zero; the gradient subtracts their mean under $q$ and multiplies by $q_i$. The differential of a loss with respect to the argument is a covector, while the residual is a vector, so an identity between the two requires a field of bilinear forms $H$ that turns the residual into a linear form, $v \mapsto H_q(X(q, x), v)$. On $\R^k$ this is the matrix $h(q)$ of Osband's principle \eqref{eq:osband-rk}, and the metric is one such field, the one for which the identity becomes $\grad_q \loss(\cdot, x) = -X(q, x)$.

\citet{ayamari2015} solve this identity by integrating the residual along the geodesics of the connection. The canonical divergence of $x$ from $q$ is
\begin{equation}\label{eq:canonical-def}
D[x:q] = \int_0^1 g\left(X(\gamma(t), x), \dot\gamma(t)\right)dt,
\end{equation}
where $\gamma$ is the geodesic from $q$ to $x$ \citep[Definition~3]{ayamari2015}; on $\R^k$ it is $\tfrac{1}{2}\|x - q\|^2$. The geodesic from $\gamma(t)$ to $x$ is the rest of $\gamma$, run in unit time, so $X(\gamma(t), x) = (1 - t)\dot\gamma(t)$ and $D[x:q] = \int_0^1 (1 - t)\,g(\dot\gamma(t), \dot\gamma(t))\,dt \ge 0$, with equality only at $x = q$ \citep[Equation~(59)]{ayamari2015}. Reversing $\gamma$ replaces the weight $1 - t$ by $t$, so $D$ is symmetric when the speed of $\gamma$ is constant, as for the Levi-Civita connection, but not in general. When the residual field is a gradient, $X(q, x) = -\grad_q f_x$ for a function $f_x$, integrating $df_x$ along $\gamma$ shows that $D[x:q] = f_x(q) - f_x(x)$, and hence $\grad_q D[x:\cdot] = -X(q, x)$ \citep[Equations~(12) and (13)]{ayamari2015}. The residual field is then called integrable. Since the simplex is contractible, this happens exactly when the lowered residual $q \mapsto \musflat{X(q, x)}$ has symmetric partial derivatives. Otherwise $D$ is still defined, but its gradient is not $-X$.

Three cases occur in this paper. For a flat connection with flat chart $\varphi$ in which the metric is the Hessian of a convex function $\psi$, $g_q(u, v) = \nabla^2\psi(\varphi(q))(d\varphi_q u, d\varphi_q v)$, the canonical divergence is the Bregman divergence in the chart, $D[x:q] = D_\psi[\varphi(x):\varphi(q)]$. Indeed, by \eqref{eq:flat-residual} its differential with respect to $q$ is $-\nabla^2\psi(\varphi(q))(\varphi(x) - \varphi(q), d\varphi_q\,\cdot) = -g_q(X(q, x), \cdot)$, the score of Section~\ref{sec:elicitability} written in the chart, and it vanishes at $q = x$. This is the case of every $\alpha$ on the cone, where $D$ is the $\alpha$-divergence of \citet[Definition~2]{ayamari2015}, and of $\alpha = \pm 1$ on the simplex, where $\psi$ is the negative entropy in the mixture chart and the log-partition function $\log(1 + \sum_{i<n} e^{\theta_i})$ in the exponential chart, so that $D[x:q] = D_{\mathrm{KL}}(x\,\|\,q)$ at $\alpha = -1$ and $D[x:q] = D_{\mathrm{KL}}(q\,\|\,x)$ at $\alpha = 1$. For the Levi-Civita connection the speed of $\gamma$ is $d(q, x)$, so $D[x:q] = \tfrac{1}{2}d(q, x)^2$ \citep[Remark~1]{ayamari2015}, and $\grad_q \tfrac{1}{2}d(\cdot, x)^2 = -\log_q(x)$ \citep[Proposition~3.10]{navasconvex}. For the other $\alpha$ the divergence $D$ is still defined by \eqref{eq:canonical-def}, but Section~\ref{sec:osband} shows that for $n \ge 3$ the residual is not integrable, and $\grad_q D[x:\cdot] \neq -X(q, x)$.

\section{The barycenter functional and its identification structure}\label{sec:barycenter}

We fix some more notation. Observations $x$ and arguments $q$ lie in the simplex, a loss $\loss(q, x)$ takes the argument first, and $\bar\loss(q, F) = \E_F\,\loss(q, X)$. For a connection with residual $X$, the expected residual is $\bar X(q, F) = \E_F X(q, X)$, a vector in the hyperplane $T_q\simplex$, which is well defined whenever $X(q, \cdot)$ is $F$-integrable.\footnote{\label{fn:integrability}On $\Ff_B$ every residual is bounded, being continuous on the compact closure of $B$. On $\Ff_0$ the mixture and Levi-Civita residuals are bounded, the latter being of length $d(q, x) < \pi$, and the exponential residual is integrable when $\E_F|\log X_i| < \infty$ for every $i$.} The definitions of strict consistency and elicitability (Definition~\ref{def:consistency}) carry over from $\R^k$ to the simplex, since they compare expected losses and use no linear structure of the space of arguments. Two classes of distributions are used. $\Ff_0$ is the class of all Borel probability measures on the open simplex. For a nonempty open set $B \subseteq \simplex$ that is geodesically convex for the connection at hand and has compact closure in the simplex, $\Ff_B$ is the class of measures with support in $B$. Both classes are convex and contain the Dirac measures at the points of $\simplex$ and of $B$, respectively.

Fix a connection $\nabla$ on the simplex with residual $X$.

\begin{definition}\label{def:barycenter}
Let $\Ff$ be a class of distributions such that, for each $F \in \Ff$, the expected residual $\bar X(q, F)$ vanishes at exactly one $q \in \simplex$. Following \citet{pennec2013}, that point is the exponential barycenter $\bary_\nabla(F)$ of $F$ for $\nabla$, and $\bary_\nabla: \Ff \to \simplex$ is the barycenter functional.
\end{definition}

On $\R^k$ the barycenter equation $\bar X(q, F) = 0$ is $\E_F X - q = 0$. On the simplex it says that the velocities at $q$ of the geodesics reaching the observations at time one average to zero. By construction, the residual is a strict $\Ff$-identification function for $\bary_\nabla$ in the sense of Section~\ref{sec:elicitability}, with values in the tangent spaces $T_q\simplex$ instead of $\R^k$. The root of $\sum_k X(\hat q, x_k) = 0$ is therefore the method-of-moments estimator of the barycenter. Where uniqueness fails on $\Ff_0$, Definition~\ref{def:barycenter} is applied to a class $\Ff_B$ on which it holds. The Levi-Civita connection of the Fisher-Rao metric requires no such restriction: its barycenter is defined on all of $\Ff_0$ (Theorem~\ref{thm:lcf-consistent}).

For the mixture connection the barycenter equation is $\E_F X - q = 0$, and $\bary_\nabla(F) = \E_F X$ is the arithmetic mean, defined on all of $\Ff_0$. For the exponential connection the expected residual $\bar X^{(e)}(q, F)_i = q_i(a_i - \sum_j q_j a_j)$, with $a_i = \E_F\log X_i - \log q_i$, vanishes exactly when all $a_i$ are equal. The barycenter is therefore the normalized geometric mean, $\bary_\nabla(F)_i \propto \exp\E_F\log X_i$, defined on the laws in $\Ff_0$ with $\E_F|\log X_i| < \infty$. For the Levi-Civita connection the barycenter equation is $\E_F\log_q(X) = 0$, the equation of \citet{karcher1977} for the Fr\'echet mean, which minimizes $q \mapsto \E_F\,d(q, X)^2$ \citep[Section~3.4]{navasconvex}; by Theorem~\ref{thm:lcf-consistent} it has exactly one solution for every $F \in \Ff_0$.

\begin{proposition}\label{prop:convex-level-sets}
The level sets of $\bary_\nabla$ are convex: if $\bary_\nabla(F_0) = \bary_\nabla(F_1) = t$ and $F_\lambda = (1 - \lambda)F_0 + \lambda F_1 \in \Ff$, then $\bary_\nabla(F_\lambda) = t$.
\end{proposition}

\begin{proof}
The expected residual is affine in the distribution, $\bar X(t, F_\lambda) = (1 - \lambda)\bar X(t, F_0) + \lambda\bar X(t, F_1) = 0$, and uniqueness in Definition~\ref{def:barycenter} implies $\bary_\nabla(F_\lambda) = t$.
\end{proof}

Every barycenter therefore passes the elicitability test: for every connection, on every convex class on which it is defined. As in Section~\ref{sec:elicitability}, the argument uses only that an expectation is linear in the distribution, however nonlinear $X(q, x)$ is in $q$ and in $x$.

The residual is not the only identification function of the barycenter; as on $\R^k$, every other one is a gauge transformation of it. The argument uses the counterpart on the simplex of the richness condition (V1) of \citet{fisslerziegel2016}, with the expected residual in the place of their identification function.

\begin{assumption}[V1]\label{ass:V1}
For every $q$ in the image of $\bary_\nabla$ there are $F_1, \dots, F_n \in \Ff$ such that $0$ is an interior point of the convex hull of $\bar X(q, F_1), \dots, \bar X(q, F_n)$ in $T_q\simplex$.
\end{assumption}

It holds on $\Ff_0$ and on every $\Ff_B$: $X(q, \cdot)$ maps a neighbourhood of $q$ diffeomorphically onto a neighbourhood of $0$ in $T_q\simplex$, so Dirac measures at suitable points near $q$ satisfy it. The following lemma is used twice, here and for Osband's principle in Section~\ref{sec:osband}. Proofs omitted from the text are in Sections~S1 and~S2 of the Supplementary Material \citep{navassupp}.

\begin{lemma}\label{lem:affine-factorization}
Let $\Ff$ be convex, let $V$ and $W$ be finite-dimensional vector spaces, and let $\Phi: \Ff \to V$ and $\Psi: \Ff \to W$ be affine maps such that $0$ is an interior point of the convex set $\Phi(\Ff)$ and $\Psi(F) = 0$ whenever $\Phi(F) = 0$. Then there is a linear map $L: V \to W$ with $\Psi = L \circ \Phi$. If moreover $\Psi(F) = 0$ only when $\Phi(F) = 0$, then $L$ is injective.
\end{lemma}

\begin{proposition}\label{prop:gauge}
Let $\Ff$ be convex and let Assumption~\ref{ass:V1} hold. For every strict $\Ff$-identification function $W$ for $\bary_\nabla$ with values in the tangent spaces there are linear isomorphisms $A_q: T_q\simplex \to T_q\simplex$ with $\bar W(q, F) = A_q\,\bar X(q, F)$ for all $q$ in the image of $\bary_\nabla$ and all $F \in \Ff$.
\end{proposition}

The proposition is Lemma~\ref{lem:affine-factorization} applied at each $q$ to $\Phi = \bar X(q, \cdot)$ and $\Psi = \bar W(q, \cdot)$, which vanish at the same distributions because both identify $\bary_\nabla$. For classes that contain the Dirac measures, such as $\Ff_0$ and $\Ff_B$, the identity at $F = \delta_x$ is the pointwise gauge transformation $W(q, x) = A_q X(q, x)$ of Section~\ref{sec:elicitability}. It is the geometric form of the characterization of identification functions by \citet*{dfzidentification}. Every moment function of the barycenter is therefore a gauge transformation of the residual, and without covariates all of them define the same estimator, since $\sum_k A_qX(q, x_k) = A_q\sum_k X(q, x_k)$.

When the residual is integrable, its canonical divergence defines a loss, $\loss(q, x) = D[x:q]$. Differentiating under the expectation, which is justified in the cases of footnote~\ref{fn:integrability} by dominated convergence, the gradient identity of Section~\ref{sec:residuals} implies
\begin{equation}\label{eq:score-identity}
d_q\bar\loss(\cdot, F) = -\musflat{\bar X(q, F)}, \quad\text{equivalently}\quad \grad_q\bar\loss(\cdot, F) = -\bar X(q, F).
\end{equation}
The differential of the expected loss is a covector and the expected residual a vector; the metric makes them comparable, by lowering the residual or, equivalently, by raising the differential to the gradient. The gradient of the expected canonical divergence is thus the negative expected residual, and the first-order condition of the M-estimator \eqref{eq:m-estimator} over the whole simplex is the barycenter equation. Over a submodel $M$ it is $P^M_q\bar X(q, F) = 0$, with $P^M_q$ the orthogonal projection onto $T_qM$ with respect to $g_q$. For the Levi-Civita connection the loss is $\tfrac{1}{2}d^2$ and the condition is $\frac{1}{N}\sum_k\log_{\hat q}(x_k) = 0$, the stationarity condition of the Levi-Civita-Fr\'echet estimator of \citet{navasconvex}. A first-order condition does not by itself render the estimator consistent; this requires that the expected loss attain its minimum at the barycenter and nowhere else.

\begin{theorem}\label{thm:lcf-consistent}
For the Levi-Civita connection of the Fisher-Rao metric, the loss $\loss(q, x) = \tfrac{1}{2}d(q, x)^2$ is strictly $\Ff_0$-consistent for the Fr\'echet mean, which is defined, unique and interior for every $F \in \Ff_0$ and is the only solution of the barycenter equation $\E_F\log_q(X) = 0$.
\end{theorem}

The proof rests on the curvature of Section~\ref{sec:residuals}. The Hessian of $\tfrac{1}{2}d(\cdot, x)^2$ at distance $d$ from $x$ has the eigenvalue $1$ along the geodesic to $x$ and $(d/2)\cot(d/2)$ across it \citep[Section~3.4]{navasconvex}, and both are positive since $d < \pi$ on the open simplex. The expected loss is therefore strictly geodesically convex on the whole simplex, which is geodesically convex. It is finite and continuous up to the boundary since $d \le \pi$, and it decreases in the directions pointing inward from the boundary, so its minimizer is interior and unique, and it is the only critical point. For the flat connections, global consistency of the Bregman losses follows from the convexity of the generator, as on $\R^k$. For the other $\alpha$-connections of the simplex the question does not arise, since by Section~\ref{sec:osband} their barycenters have no consistent loss in the regularity class of Osband's principle.

\section{Geometric approach to Osband's principle}\label{sec:osband}

On $\R^k$ the matrix $h(q)$ of Osband's principle \eqref{eq:osband-rk} is determined by the loss. On the simplex it becomes a tensor field $H$, with a bilinear form $H_q: T_q\simplex \times T_q\simplex \to \R$ at each $q$, turning the residual into a linear form as in Section~\ref{sec:residuals}. Let $\Ff$ be convex and let $A$ be the interior of the image of $\bary_\nabla$.

\begin{theorem}\label{thm:osband}
Let Assumption~\ref{ass:V1} hold, and let $\loss$ be a strictly $\Ff$-consistent loss for $\bary_\nabla$ such that $q \mapsto \bar\loss(q, F)$ is continuously differentiable on $A$ for every $F \in \Ff$. Then there is a unique tensor field $H$ of bilinear forms on the tangent spaces $T_q\simplex$, $q \in A$, with
\begin{equation}\label{eq:osband-simplex}
d_q\bar\loss(\cdot, F) = -H_q(\bar X(q, F), \cdot) \quad\text{for all } q \in A \text{ and } F \in \Ff.
\end{equation}
\end{theorem}

The condition on $\bar\loss$ is the counterpart of Assumption (S1) of \citet{fisslerziegel2016}, and Theorem~\ref{thm:osband} is the geometric form of their Theorem~3.2. Its proof applies Lemma~\ref{lem:affine-factorization} at each $q \in A$ to two affine maps: $\Phi = \bar X(q, \cdot)$ from $\Ff$ to $T_q\simplex$, and $\Psi = d_q\bar\loss(\cdot, \cdot)$ from $\Ff$ to its dual. By Assumption~\ref{ass:V1} the zero vector is an interior point of $\Phi(\Ff)$. If $\bar X(q, F) = 0$, then $q = \bary_\nabla(F)$, the unique minimizer of $\bar\loss(\cdot, F)$ by strict consistency, so $d_q\bar\loss(\cdot, F) = 0$. The lemma produces a linear map $L_q: T_q\simplex \to T_q^{\ast}\simplex$ with $\Psi = L_q \circ \Phi$, and we set $H_q(u, v) = -(L_qu)(v)$. The expected residuals span $T_q\simplex$, so $H_q$ is unique. By the second part of the lemma, $L_q$ would be injective, and $H_q$ nondegenerate, if $d_q\bar\loss(\cdot, F)$ vanished only for $q = \bary_\nabla(F)$. Strict consistency makes the barycenter the unique minimizer of the expected loss but does not exclude other critical points, so $H_q$ may be degenerate. We call $H$ the \emph{Osband tensor} of $\loss$.

The definition of $H$ uses the loss and the connection, through its residual, but it does not use the metric. Take a flat connection with flat chart $\varphi$. By \eqref{eq:flat-residual} its expected residual vanishes exactly when $\E_F\varphi(X) = \varphi(q)$, so its barycenter is the mean in the chart, $\bary_\nabla(F) = \varphi^{-1}(\E_F\varphi(X))$. For a strictly convex $\psi$ on the convex image of $\varphi$, the Bregman loss $\loss(q, x) = D_\psi[\varphi(x):\varphi(q)]$ is therefore strictly consistent for it by Section~\ref{sec:elicitability}. Differentiating as at the end of Section~\ref{sec:residuals} and taking expectations, $d_q\bar\loss(\cdot, F) = -\nabla^2\psi(\varphi(q))(d\varphi_q\bar X(q, F), d\varphi_q\,\cdot)$. The Osband tensor of the Bregman loss is therefore the Hessian of the generator in the flat chart, $H_q(u, v) = \nabla^2\psi(\varphi(q))(d\varphi_qu, d\varphi_qv)$. On $\R^k$, with $\varphi$ the identity, this is the matrix $h = \nabla^2\psi$ of Section~\ref{sec:elicitability}. When the canonical divergence of a connection is a consistent loss, \eqref{eq:score-identity} shows that its Osband tensor is the metric, $H = g$.

A metric $g$ identifies each bilinear form $H_q$ with the endomorphism $h_q: T_q\simplex \to T_q\simplex$, $u \mapsto \mussharp{(H_q(u, \cdot))}$, so that $g_q(h_qu, v) = H_q(u, v)$. Raising \eqref{eq:osband-simplex} yields
\begin{equation}\label{eq:osband-grad}
\grad_q\bar\loss(\cdot, F) = -h_q\bar X(q, F),
\end{equation}
which is Osband's principle \eqref{eq:osband-rk} with $h_q$ in the place of the matrix $h(q)$; for the Euclidean metric on $\R^k$ the two coincide. The endomorphism depends on the metric; the Osband tensor does not. Where $H_q$ is nondegenerate, $h_q$ is invertible, and \eqref{eq:osband-grad} says that the negative gradient of the expected loss is the expected residual after the gauge transformation $h_q$ of Section~\ref{sec:elicitability}. Choosing a consistent loss with nondegenerate Osband tensor is thus choosing a gauge of the residual, and the canonical divergence, with $H = g$, has the identity gauge, as $\tfrac{1}{2}\|x - q\|^2$ does on $\R^k$. Over the whole simplex an invertible $h_q$ cancels from the first-order condition, leaving the barycenter equation. Over a submodel $M$ the first-order condition $P^M_q h_q\bar X(q, F) = 0$ depends on the gauge: all these equations are solved at the barycenter when it lies in $M$, and their solutions differ otherwise (Section~\ref{sec:regression}). Not every field of endomorphisms is the gauge of a loss: the one-form $q \mapsto H_q(\bar X(q, F), \cdot)$ must be exact, with potential $-\bar\loss(\cdot, F)$, and on a contractible $A$, such as the simplex, this holds exactly when it is closed.

The classes $\Ff_0$ and $\Ff_B$ of Section~\ref{sec:barycenter} contain the Dirac measures. So the counterpart of Assumption (F1) of \citet{fisslerziegel2016}, approximation of Dirac measures within the class, holds trivially. At $F = \delta_x$, \eqref{eq:osband-simplex} is an identity for the loss itself,
\begin{equation}\label{eq:osband-pointwise}
d_q\loss(\cdot, x) = -H_q(X(q, x), \cdot), \quad q \in A.
\end{equation}
We say that $\loss$ belongs to the \emph{regularity class of Osband's principle} if it satisfies the hypotheses of Theorem~\ref{thm:osband} and $\loss(\cdot, x)$ is twice continuously differentiable on $A$ for every $x$.

\begin{proposition}\label{prop:closedness}
Let $\loss$ belong to the regularity class of Osband's principle. Then its Osband tensor $H$ is continuously differentiable, and for every $x$ the one-form $\omega_x: q \mapsto H_q(X(q, x), \cdot)$ on $A$ is closed. Conversely, let $A$ be contractible and let $H$ be a continuously differentiable tensor field of bilinear forms on $A$ such that every $\omega_x$ is closed. Then $\loss(q, x) = -\int_x^q \omega_x$, along any path in $A$ from $x$ to $q$, is well defined, vanishes at $q = x$ and satisfies \eqref{eq:osband-pointwise}.
\end{proposition}

\begin{proof}
Near $q_0 \in A$, choose observations $x_1, \dots, x_{n-1}$ whose residuals $X(q_0, x_i)$ form a basis of $T_{q_0}\simplex$. They remain a basis near $q_0$, and by \eqref{eq:osband-pointwise}, $H_q$ solves the linear system $H_q(X(q, x_i), v) = -d_q\loss(\cdot, x_i)(v)$, whose coefficients are continuously differentiable in $q$. By \eqref{eq:osband-pointwise}, $\omega_x = -d\loss(\cdot, x)$ is exact, hence closed. Conversely, closed forms on a contractible set are exact by the Poincar\'e lemma. The path integral is therefore independent of the path, and its differential with respect to $q$ is $-\omega_x$, which is \eqref{eq:osband-pointwise}.
\end{proof}

The expected differential of the loss in the converse is $-H_q(\bar X(q, F), \cdot)$, which vanishes at the barycenter. Consistency requires in addition that the expected loss attain its minimum there and nowhere else. Convexity of the generator ensures this for the Bregman losses, and Theorem~\ref{thm:lcf-consistent} ensures it for $\tfrac{1}{2}d^2$. Closedness of the forms $\omega_x$ is thus a condition on the pair $(H, X)$ in which the loss no longer appears. Every loss in the regularity class is recovered from its Osband tensor by integration, up to the additive function $\loss(x, x)$.

For $H = g$, \eqref{eq:osband-pointwise} becomes $d_q\loss(\cdot, x) = -\musflat{X(q, x)}$. The residual is then integrable with potential $\loss(\cdot, x)$, and by Section~\ref{sec:residuals}, $\loss(q, x) - \loss(x, x) = D[x:q]$.

\begin{proposition}\label{prop:canonical}
A loss in the regularity class of Osband's principle has Osband tensor $g$ if and only if the residual is integrable on $A$ and $\loss(q, x) = D[x:q] + a(x)$, with $a(x) = \loss(x, x)$. The canonical divergence belongs to this class only if the residual is integrable.
\end{proposition}

The converse direction of the equivalence is the gradient identity of Section~\ref{sec:residuals}. For the last statement, suppose $D$ is in the class. By \eqref{eq:canonical-def}, $D[x:q] = \tfrac{1}{2}g_q(X(q, x), X(q, x)) + O(\|X(q, x)\|^3)$, so $d_qD[x:\cdot] = -g_q(X(q, x), \cdot) + O(\|X(q, x)\|^2)$. Comparing with \eqref{eq:osband-pointwise}, which is linear in $X(q, x)$, yields $H_q = g_q$, and the first statement then makes the residual integrable. \textbf{A canonical divergence that is a consistent loss has the metric as its Osband tensor, and the residual of its connection is integrable.}

Ay and Amari's integral \eqref{eq:canonical-def} is thus the reconstruction of a loss from the Osband tensor $g$, along the geodesics of the connection and based at the observation. The normalization $D[x:x] = 0$ fixes the additive function $a$, the same freedom as in the Bregman losses \eqref{eq:bregman-def}. Up to $a$, the canonical divergence is the only loss whose gauge is the identity. In the flat cases the consistent losses are determined only up to the choice of a convex generator, and the canonical divergence is the Bregman loss whose Hessian is the metric.

On $\R^k$ the residual $x - q$ has derivative $-u$ in every direction $u$, and the forms $\omega_x$ are closed exactly when $H$ is the Hessian of a function (Section~\ref{sec:elicitability}). For a curved connection the derivative of the residual contains the curvature. The $\alpha$-connections are torsion-free, since their Christoffel symbols are symmetric. So $d\omega(u, v) = (\nabla_u\omega)(v) - (\nabla_v\omega)(u)$ for every one-form $\omega$, and by the Leibniz rule
\begin{equation}\label{eq:closedness-expanded}
d\omega_x(u, v) = (\nabla_uH)(X, v) - (\nabla_vH)(X, u) + H(\nabla_uX, v) - H(\nabla_vX, u), \quad X = X(\cdot, x).
\end{equation}
To compute $\nabla_uX(\cdot, x)$, move the argument along a curve $q(s)$ with $q(0) = q$ and $\dot q(0) = u$, and keep the observation $x$ fixed. Each $q(s)$ is joined to $x$ by a geodesic. The rate at which this family of geodesics changes with $s$ is a vector field $J$ along the geodesic $\gamma$ from $q$ to $x$, called a \emph{Jacobi field}. It equals $u$ at $q$ and vanishes at $x$, where all the geodesics end. It solves the Jacobi equation $\nabla_{\dot\gamma}\nabla_{\dot\gamma}J + R(J, \dot\gamma)\dot\gamma = 0$, through which the curvature $R$ of the connection enters. Its derivative at $q$ is $\nabla_uX(\cdot, x)$. Expanding $J$ in powers of the time along $\gamma$ and imposing $J = 0$ at $x$ yields the following well-known expansion; for torsion-free connections it follows from the expansion of the inverse exponential map at a moving base point in \citet[Theorem~3]{pennec2019}, differentiated with respect to the base point.

\begin{lemma}\label{lem:jacobi}
Let $\nabla$ be a torsion-free connection with curvature tensor $R$ and residual $X$. For $q \in A$, $u \in T_q\simplex$ and $x$ in a neighbourhood of $q$ on which the residual is defined, with $L = X(q, x)$,
\begin{equation}\label{eq:jacobi-expansion}
\nabla_uX(\cdot, x)|_q = -u + \tfrac{1}{3}R(u, L)L + \tfrac{1}{12}(\nabla_LR)(u, L)L + O(\|L\|^4).
\end{equation}
\end{lemma}

For the Levi-Civita connection of the Fisher-Rao metric, \eqref{eq:jacobi-expansion} can be compared with Section~\ref{sec:barycenter}. At $\alpha = 0$, \eqref{eq:alpha-curvature} becomes $\bar R^{(0)}(u, L)L = \tfrac{1}{4}(g(L, L)u - g(u, L)L)$. This curvature is constant, so its covariant derivative vanishes, and \eqref{eq:jacobi-expansion} becomes $\nabla_uX = -u + \tfrac{1}{12}(g(L, L)u - g(u, L)L) + O(\|L\|^4)$. Along $L$ the correction vanishes and $\nabla_LX = -L$, while across $L$, for $g(u, L) = 0$, one obtains $\nabla_uX = -(1 - \|L\|^2/12)u$. On the other hand, $X(\cdot, x) = \log_{\cdot}(x) = -\grad\tfrac{1}{2}d(\cdot, x)^2$ by Section~\ref{sec:residuals}, so $\nabla_uX = -\Hess\tfrac{1}{2}d(\cdot, x)^2\,u$. The eigenvalues of this Hessian are $1$ along $L$ and $(d/2)\cot(d/2)$ across it, with $d = d(q, x) = \|L\|$, and since $(d/2)\cot(d/2) = 1 - d^2/12 + O(d^4)$, the two descriptions agree.

\begin{theorem}\label{thm:rigidity}
Let $n \ge 3$ and $\alpha^2 \ne 1$, and let $H$ be a continuously differentiable tensor field of bilinear forms on a connected open set $A \subseteq \simplex$ such that, for every $q \in A$ and every $x$ in a neighbourhood of $q$, the one-form $\omega_x = H(\bar X^{(\alpha)}(\cdot, x), \cdot)$ is closed near $q$. Then $H = c\,g$ on $A$ for a constant $c \in \R$.
\end{theorem}

\begin{proof}
Fix $q \in A$ and $u, v \in T_q\simplex$. As $x$ runs through a neighbourhood of $q$, $L = \bar X^{(\alpha)}(q, x)$ runs through a neighbourhood $U$ of $0$, and by \eqref{eq:jacobi-expansion} the right-hand side of \eqref{eq:closedness-expanded} is $E_0 + E_1(L) + E_2(L) + O(\|L\|^3)$, with $E_0 = H(v, u) - H(u, v)$, $E_1(L) = (\nabla_uH)(L, v) - (\nabla_vH)(L, u)$ and $E_2(L) = \tfrac{1}{3}(H(\bar R^{(\alpha)}(u, L)L, v) - H(\bar R^{(\alpha)}(v, L)L, u))$, homogeneous of degrees $0$, $1$ and $2$ in $L$. Closedness along $tL$ for small $t$ makes each of them vanish. $E_0 = 0$ says that $H$ is symmetric. By \eqref{eq:alpha-curvature}, $E_2(L) = \frac{1 - \alpha^2}{12}(g(L, v)H(L, u) - g(L, u)H(L, v))$, and its vanishing at $u = L$ yields $H(L, \cdot) = \lambda(L)\,g(L, \cdot)$ with $\lambda(L) = H(L, L)/g(L, L)$. Every $L \in U \setminus \{0\}$ is thus an eigenvector of the endomorphism $h_q$ of \eqref{eq:osband-grad}; a finite union of proper subspaces has empty interior, so one eigenspace is all of $T_q\simplex$ and $H_q = \lambda(q)\,g_q$. Finally, $(\bar\nabla^{(\alpha)}_ug)(v, w) = \alpha\,T(u, v, w)$ with the Amari-Chentsov tensor $T(u, v, w) = \sum_i u_iv_iw_i/q_i^2$, which is symmetric in its three arguments \citep{amari2000}, so $E_1(L) = d\lambda(u)\,g(L, v) - d\lambda(v)\,g(L, u)$. For $v = L$ and $u$ orthogonal to $L$, which exists since $\dim T_q\simplex = n - 1 \ge 2$, this is $d\lambda(u)\,g(L, L)$. Hence $d\lambda = 0$, and $\lambda$ is constant on the connected set $A$.
\end{proof}

The condition $n \ge 3$ is necessary: for $n = 2$ the simplex is one-dimensional, every one-form on it is closed, and every $H$ passes the closedness test.

At $\alpha = \pm 1$ the curvature vanishes and, by \eqref{eq:flat-residual}, $\nabla_uX = -u$ exactly in the flat chart. The right-hand side of \eqref{eq:closedness-expanded} then vanishes for all $x$ exactly when $H$ is symmetric and $\nabla H$ is symmetric in its three arguments. In the flat chart, where the Christoffel symbols vanish, this makes $H$ the Hessian of a function $\psi$, and the consistent losses are the Bregman losses of Section~\ref{sec:elicitability} written in the flat chart.

By Theorem~\ref{thm:rigidity} the only possibility left for $\alpha^2 \ne 1$ is $H = g$, and by Proposition~\ref{prop:canonical} it belongs to a consistent loss only if the residual is integrable. At $\alpha = 0$ the residual is integrable, with potential $\tfrac{1}{2}d(\cdot, x)^2$ (Section~\ref{sec:residuals}). For $\alpha \notin \{-1, 0, 1\}$ it is not.

\begin{proposition}\label{prop:nonintegrable}
Let $n \ge 3$ and $\alpha \notin \{-1, 0, 1\}$. For $q \in A$, $u, v, L \in T_q\simplex$ and $x = \exp_q(tL)$, write $r_i = L_i/q_i$. Then
\begin{equation}\label{eq:cubic-obstruction}
d\musflat{\bar X^{(\alpha)}(\cdot, x)}(u, v) = t^3\,\frac{(1 - \alpha^2)\alpha}{48}\,B_L(u, v) + O(t^4),
\end{equation}
with $B_L(u, v) = T(L, v, L)\,g(L, u) - T(L, u, L)\,g(L, v) = \sum_i r_i^2v_i\sum_j r_ju_j - \sum_i r_i^2u_i\sum_j r_jv_j$.
\end{proposition}

The bracket $B_L$ vanishes for all $u, v \in T_q\simplex$ exactly when the linear forms $u \mapsto \sum_i r_i^2u_i$ and $u \mapsto \sum_i r_iu_i$ are linearly dependent on the tangent space. Since tangent vectors satisfy $\sum_i u_i = 0$, this means $r_i^2 = \lambda r_i + \mu$ for all $i$, with constants $\lambda$ and $\mu$, so the ratios $r_i$ take at most two distinct values. For $n = 2$ there are only two ratios, and $B_L$ always vanishes, in line with the closedness of every one-form on the one-dimensional simplex. For $n \ge 3$ three of the ratios can be distinct, and then $B_L \ne 0$. The residual $\bar X^{(\alpha)}$ is therefore not integrable.

The expansion sets $H = g$ in \eqref{eq:closedness-expanded} and inserts \eqref{eq:jacobi-expansion}. The terms with $\bar\nabla^{(\alpha)}g$ are $\alpha(T(u, L, v) - T(v, L, u))$, and they cancel because $T$ is symmetric in its three arguments. The term of degree zero is $g(-u, v) - g(-v, u) = 0$. The term of degree two is $\tfrac{1}{3}(g(\bar R^{(\alpha)}(u, L)L, v) - g(\bar R^{(\alpha)}(v, L)L, u))$. By \eqref{eq:alpha-curvature}, $g(\bar R^{(\alpha)}(u, L)L, v) = \frac{1 - \alpha^2}{4}(g(L, L)\,g(u, v) - g(u, L)\,g(L, v))$, which is symmetric in $u$ and $v$, so the difference vanishes. The first term that remains has degree three, $\tfrac{1}{12}(g((\nabla_L\bar R^{(\alpha)})(u, L)L, v) - g((\nabla_L\bar R^{(\alpha)})(v, L)L, u))$. Differentiating \eqref{eq:alpha-curvature} with $\bar\nabla^{(\alpha)}g = \alpha T$ yields
\[
(\nabla_L\bar R^{(\alpha)})(u, L)L = \frac{(1 - \alpha^2)\alpha}{4}(T(L, L, L)\,u - T(L, u, L)\,L),
\]
and antisymmetrizing in $u$ and $v$ yields \eqref{eq:cubic-obstruction}. Numerical checks at $\alpha = \pm 1/2$ are in Section~S3 of the Supplementary Material \citep{navassupp}.

\begin{corollary}\label{cor:regimes}
Let $n \ge 3$, and let $\loss$ be a loss in the regularity class of Osband's principle for the barycenter of $\bar\nabla^{(\alpha)}$ on $\Ff_0$ or on a class $\Ff_B$. Then, up to an additive function of the observation,
\begin{enumerate}[label=(\roman*)]
\item at $\alpha = \pm 1$, $\loss(q, x) = D_\psi[\varphi(x):\varphi(q)]$ for a strictly convex $\psi$ in the flat chart $\varphi$;
\item at $\alpha = 0$, $\loss(q, x) = \tfrac{c}{2}d(q, x)^2$ for a constant $c > 0$;
\item for $\alpha \notin \{-1, 0, 1\}$, no such loss exists.
\end{enumerate}
\end{corollary}

For $\alpha^2 \ne 1$, Proposition~\ref{prop:closedness} and Theorem~\ref{thm:rigidity} yield $H = c\,g$. By \eqref{eq:osband-pointwise}, the Hessian of $\loss(\cdot, x)$ at its minimum $q = x$ is $c\,g_x$, so $c \ge 0$; and $c = 0$ would make $\loss(\cdot, x)$ constant, so $c > 0$. Proposition~\ref{prop:canonical} applied to $\loss/c$ then requires an integrable residual, which Proposition~\ref{prop:nonintegrable} excludes for $\alpha \notin \{-1, 0, 1\}$. At $\alpha = 0$ it yields $\loss = c\,D + a$, which is consistent by Theorem~\ref{thm:lcf-consistent}. In case (iii) the barycenter equation $\sum_k \bar X^{(\alpha)}(\hat q, x_k) = 0$ does not come from the first-order conditions of a minimization. It remains a method-of-moments estimator, with the limit law of Section~\ref{sec:regression}. Whether a loss outside the regularity class of Osband's principle is strictly consistent for this barycenter, or whether the barycenter is not elicitable at all, is open.

Proposition~\ref{prop:nonintegrable} corrects the statement of \citet[Section~6]{ayamari2015} that the restriction $\bar\nabla^{(\alpha)}$ of the $\alpha$-connection to the simplex ``is still integrable'', although the simplex is not dually flat for $\alpha \notin \{-1, 1\}$. The statement holds at $\alpha = 0$ and fails for every $\alpha \notin \{-1, 0, 1\}$ when $n \ge 3$. The reason lies in the restriction factor of their Equation~(36), discussed in Section~\ref{sec:alpha}.

\section{Extrinsic and intrinsic estimation along the \texorpdfstring{$\alpha$}{alpha}-family}\label{sec:alpha}

For $\alpha^2 \neq 1$, i.e. where $\bar\nabla^{(\alpha)}$ is curved, the residual of the simplex differs from the projected residual of the cone by a factor. Normalizing the cone geodesic from $q$ to $x$ produces a curve $\sigma_{q,x}$ in the simplex from $q$ to $x$. \citet[Section~3.3]{ayamari2015} show that the geodesic of $\bar\nabla^{(\alpha)}$ from $q$ to $x$ is this curve reparametrized in time, $t \mapsto \sigma_{q,x}(\tau_{q,x}(t))$ with $\tau_{q,x}(0) = 0$ and $\tau_{q,x}(1) = 1$. In the chart $p \mapsto ((p_1/p_n)^{(1-\alpha)/2}, \dots, (p_{n-1}/p_n)^{(1-\alpha)/2})$ the curves $\sigma_{q,x}$ are straight segments in the open positive orthant of $\R^{n-1}$. So $\bar\nabla^{(\alpha)}$ is projectively flat, meaning its geodesics are straight lines of a chart up to parametrization, and any two points are joined by exactly one geodesic. This is an instance of a general fact. With \eqref{eq:alpha-curvature}, $(\simplex, g, \bar\nabla^{(\alpha)})$ is a statistical manifold of constant curvature $(1 - \alpha^2)/4$ in the sense of \citet{kurose1990}, and such manifolds are projectively flat \citep[Theorem~2.6]{kobayashiohno2022}. Normalization differentiates to $P_q$ at a point of the simplex, so the chain rule yields their Equation~(36),
\begin{equation}\label{eq:invexp-restricted}
\bar X^{(\alpha)}(q, x) = \dot\tau_{q,x}(0)\,P_qX^{(\alpha)}(q, x).
\end{equation}
The factor $\dot\tau_{q,x}(0)$ depends on both points. It equals $1$ at $\alpha = \pm 1$, where the normalized geodesics need no reparametrization, in agreement with $\bar X = P_qX$ in Section~\ref{sec:residuals}. At $\alpha = 0$ the residual is the Riemannian logarithm $\log_q(x)$, of length $d(q, x)$. Mapping the great-circle arc from $\Pi(q)$ to $\Pi(x)$ back to the simplex yields, with $d = d(q, x)$,
\[
\log_q(x)_i = \frac{d}{\sin(d/2)}\left(\sqrt{q_ix_i} - q_i\sum_j \sqrt{q_jx_j}\right),
\]
so that $\dot\tau_{q,x}(0) = (d/2)/\sin(d/2)$, which exceeds $1$ for $x \neq q$. A factor that depends on the observation is not a gauge transformation in the sense of Section~\ref{sec:elicitability}, as it can change where the expected residual vanishes.

The cone is flat for $\nabla^{(\alpha)}$, so its residual is integrable, with the $\alpha$-divergence of \citet{ayamari2015} as potential (Section~\ref{sec:residuals}). By \eqref{eq:grad-simplex} its projection is the negative gradient on the simplex of the restricted divergence,
\[
P_qX^{(\alpha)}(q, x) = -\grad_q D^{(\alpha)}[x:\cdot], \quad D^{(\alpha)}[x:q] = \frac{4}{1 - \alpha^2}\left(1 - \sum_i x_i^{(1-\alpha)/2}q_i^{(1+\alpha)/2}\right).
\]
By \eqref{eq:invexp-restricted} the residual of the simplex is proportional to a gradient,
\begin{equation}\label{eq:tau-factor}
\bar X^{(\alpha)}(q, x) = -\dot\tau_{q,x}(0)\,\grad_q D^{(\alpha)}[x:\cdot].
\end{equation}
This proportionality is all that the geodesic projection property of Ay and Amari (their Proposition~6) requires, since there the factor may depend on both points, and it holds for every $\alpha$. Integrability requires more. The one-form $-\dot\tau_{\cdot,x}(0)\,dD^{(\alpha)}[x:\cdot]$ is closed exactly when $d\dot\tau_{\cdot,x}(0) \wedge dD^{(\alpha)}[x:\cdot] = 0$, that is, when $\dot\tau_{\cdot,x}(0)$ is locally a function of $D^{(\alpha)}[x:\cdot]$. At $\alpha = 0$ this holds: $D^{(0)}[x:q] = 2\sum_i(\sqrt{x_i} - \sqrt{q_i})^2$, twice the squared Hellinger distance, and $\dot\tau_{q,x}(0) = (d/2)/\sin(d/2)$ with $d = 2\arccos\sum_i\sqrt{q_ix_i}$ both depend on $q$ only through $\sum_i\sqrt{q_ix_i}$. For $\alpha \notin \{-1, 0, 1\}$ and $n \ge 3$ it fails, by Proposition~\ref{prop:nonintegrable}. What restriction to the simplex preserves is thus the integrability of the projected cone residual, whose potential is the restricted $\alpha$-divergence, together with the proportionality \eqref{eq:tau-factor} that Proposition~6 of Ay and Amari uses. The factor $\dot\tau_{q,x}(0)$ separates the projected cone residual from the residual of $\bar\nabla^{(\alpha)}$.

At the ends of the family of restricted $\alpha$-divergences, the Kullback-Leibler divergences are $D^{(-1)}[x:q] = D_{\mathrm{KL}}(x\,\|\,q)$ and $D^{(1)}[x:q] = D_{\mathrm{KL}}(q\,\|\,x)$ (Section~\ref{sec:residuals}). They elicit the barycenters of Section~\ref{sec:barycenter}: the arithmetic mean $\E_F X$, and the normalized geometric mean, whose components are proportional to $\exp\E_F\log X_i$. Both are normalized power means $(\E_F X_i^r)^{1/r}$, with exponent $r = 1$ and, as a limit, $r = 0$. Between the ends the exponent $r = (1 - \alpha)/2$ runs through $(0, 1)$, and the restricted $\alpha$-divergence, as a loss $\loss(q, x) = D^{(\alpha)}[x:q]$, elicits the corresponding power mean. For $F \in \Ff_0$ let $m_i(F) = \E_F X_i^{(1-\alpha)/2}$, which lies in $(0, 1]$.

\begin{proposition}\label{prop:power-mean}
For $\alpha \in (-1, 1)$, the loss $D^{(\alpha)}[x:q]$ is strictly $\Ff_0$-consistent for the normalized power mean
\begin{equation}\label{eq:power-mean}
\Pi^{(\alpha)}(F)_i = \frac{m_i(F)^{2/(1-\alpha)}}{\sum_j m_j(F)^{2/(1-\alpha)}},
\end{equation}
and its M-estimator is the plug-in power mean $\hat q_i \propto (\frac{1}{N}\sum_k x_{k,i}^{(1-\alpha)/2})^{2/(1-\alpha)}$.
\end{proposition}

The expected loss is $\frac{4}{1 - \alpha^2}(1 - \sum_i m_i(F)\,q_i^{(1+\alpha)/2})$, so it depends on $F$ only through $m(F)$. Minimizing it over the simplex means maximizing $\sum_i m_iq_i^{(1+\alpha)/2}$, which is strictly concave since $(1 + \alpha)/2 \in (0, 1)$. The partial derivatives of this sum diverge at the boundary, so the maximizer is interior, and the Lagrange condition $m_iq_i^{(\alpha-1)/2} = \mathrm{const}$ yields \eqref{eq:power-mean}. The plug-in power mean is a method-of-moments estimator in the original sense: it replaces the moments $m_i(F)$ by their sample averages and transforms them.

The power mean is the barycenter of the flat connection $\nabla^{(\alpha)}$ of the cone, normalized to the simplex. In the flat chart $m \mapsto m^{(1-\alpha)/2}$ of Section~\ref{sec:residuals}, the barycenter of $F$ on the cone is the mean in the chart, which has coordinates $m_i(F)$, as described for flat connections in Section~\ref{sec:osband}. Normalizing the point with these coordinates yields \eqref{eq:power-mean}.

By \eqref{eq:tau-factor}, the first-order condition of this M-estimator over the simplex is $\sum_k \grad_q D^{(\alpha)}[x_k:\cdot] = 0$, that is,
\begin{equation}\label{eq:weighted-score}
\sum_k \dot\tau_{\hat q,x_k}(0)^{-1}\,\bar X^{(\alpha)}(\hat q, x_k) = 0.
\end{equation}
This is a weighted barycenter equation. The barycenter of $\bar\nabla^{(\alpha)}$, the intrinsic functional, solves the unweighted equation $\sum_k \bar X^{(\alpha)}(\hat q, x_k) = 0$. The two coincide at $\alpha = \pm 1$, where $\dot\tau \equiv 1$. Otherwise the weights are not constant, and the two equations have different solutions in general.

At $\alpha = 0$ the weights $\sin(d_k/2)/(d_k/2)$ decrease with the distance $d_k = d(\hat q, x_k)$, so the extrinsic equation discounts distant observations relative to the intrinsic one. Its solution is the power mean $\Pi^{(0)}(F)_i \propto (\E_F\sqrt{X_i})^2$, elicited by twice the squared Hellinger distance. The intrinsic solution is the Fr\'echet mean, elicited by $\tfrac{1}{2}d^2$ (Theorem~\ref{thm:lcf-consistent}). Under the isometry $\Pi$ these are the extrinsic and intrinsic means of \citet{bhattacharya2003}. The first projects the Euclidean mean of the points $2\sqrt{X}$ onto the sphere; the second minimizes the expected squared geodesic distance on it. Both are elicitable, but not by the same loss. A forecast of the Fr\'echet mean evaluated with twice the squared Hellinger distance, or a forecast of the power mean evaluated with $\tfrac{1}{2}d^2$, is in general scored against the wrong target.

For $\alpha \notin \{-1, 0, 1\}$ the split is sharper. The extrinsic power mean keeps its consistent loss, while the intrinsic barycenter has none in the regularity class of Osband's principle, by Corollary~\ref{cor:regimes}.

\section{Regression on the simplex: model consistency and instruments}\label{sec:regression}

In the simplest setting a real response $Y$ is explained by covariates $Z \in \R^r$ through the linear model $m(z, \theta) = z\T\theta$, $\theta \in \R^r$. The model is correctly specified if $\E[Y \,|\, Z] = Z\T\theta_0$. Least squares minimizes $\frac{1}{N}\sum_{k=1}^N (y_k - z_k\T\theta)^2$, and its first-order condition, the normal equations $\sum_k z_k(y_k - z_k\T\hat\theta) = 0$, multiplies the residual of each observation by its covariates. It estimates $\theta_0$ because the squared loss is consistent for the mean. At each value of $Z$ the conditional expected loss is minimized by the conditional mean, and the model attains the conditional mean at $\theta_0$.

In this example the conditional mean is the \emph{target}, the squared loss is a \emph{loss} consistent for it, and the linear function $z\T\theta$ is the \emph{model}. The target determines which losses are consistent for it, while the model is chosen separately. For a binary response the target $\E[Y \,|\, Z] = P(Y = 1 \,|\, Z)$ stays the same, while logit, $p = e^{z\T\theta}/(1 + e^{z\T\theta})$, probit, $p = \Phi(z\T\theta)$, and the linear probability model, $p = z\T\theta$, are three models for it. Each of them can be fitted with the Bernoulli log-likelihood or with the squared loss, both consistent for the mean.

On the simplex the target is the barycenter of a connection, the losses are those of Sections~\ref{sec:barycenter} to~\ref{sec:alpha}, and the model is a smooth map $m: \R^r \times \Theta \to \simplex$, $\Theta \subseteq \R^s$. A connection produces a natural class of models, the geodesic models
\begin{equation}\label{eq:geodesic-model}
m(z, \theta) = \exp_p(z_1v_1 + \dots + z_rv_r), \quad \theta = (p, v_1, \dots, v_r), \quad v_j \in T_p\simplex,
\end{equation}
with $s = (r + 1)(n - 1)$ parameters in coordinates. \citet{fletcher2013} introduced them for the Levi-Civita connection of a Riemannian manifold and one covariate, and on the Fisher-Rao simplex their fitted values follow great circles under $\Pi$. For the mixture connection they are the linear models $p + \sum_j z_jv_j$, for covariates at which these stay in the simplex. For the exponential connection they are the log-linear models $m(z, \theta)_i \propto p_i\exp(\sum_j z_jv_{j,i}/p_i)$. These are the multinomial logistic models, with intercepts $\log p_i$ and slopes $v_{j,i}/p_i$, the standard way of turning categories into compositions \citep{mcfadden1974,fernandesmachado2025}. For $n = 2$ the mixture and exponential geodesic models are the linear probability model and the logit. The probit is the geodesic model of the connection that is flat in the chart $p \mapsto \Phi^{-1}(p)$. The connection of the model need not be that of the target.

A model on the simplex is correctly specified for a functional $\Gamma$ if there is a unique $\theta_0$ in the interior of $\Theta$ with $m(Z, \theta_0) = \Gamma(F_{Y|Z})$ almost surely. A loss is \emph{strictly conditionally model consistent} if $\theta_0$ is almost surely the unique minimizer of $\theta \mapsto \E[\loss(m(Z, \theta), Y) \,|\, Z]$, and \emph{strictly unconditionally model consistent} if $\theta_0$ is the unique minimizer of the unconditional expectation. \citetalias{dimitriadis2024} prove in their Theorem~1 that a strictly consistent loss for $\Gamma$ is strictly conditionally model consistent, and that conditional model consistency implies unconditional model consistency. The converse implications hold under their Assumptions~2 and~3, the first of them for a modified class of conditional laws. These assumptions ask for a surjective model, a conditional expected loss continuous in the parameter, and a rich class of joint laws. Their functional takes values in an arbitrary subset of a Euclidean space, and the proof uses no further structure of that subset, so the theorem applies with the simplex as that subset.

Categorical responses lie on the boundary of the simplex, while Sections~\ref{sec:residuals} to~\ref{sec:alpha} are set on the open simplex. Two of the connections extend to them without change. For the mixture connection, the residual $y - q$ and the loss $D_{\mathrm{KL}}(y\,\|\,q)$, with $0\log 0 = 0$, are defined for $y$ in the closed simplex and $q$ in the open one. For the Levi-Civita connection, the vertices lie at distance $d(q, e_c) = 2\arccos\sqrt{q_c} < \pi$, and $\tfrac{1}{2}d^2$ and $\log_q$ extend continuously to the closed simplex \citep[Proposition~3.11 and Section~3.4]{navasconvex}. On the laws on the closed simplex under which every component is positive with positive probability, the barycenters of both connections lie in the open simplex, and $D_{\mathrm{KL}}(y\,\|\,q)$ and $\tfrac{1}{2}d^2$ remain strictly consistent for them. The exponential residual contains $\log(y_i/q_i)$ and does not extend.

By Theorem~1 of \citetalias{dimitriadis2024}, parts (i) and (iii), the consistent losses of Sections~\ref{sec:barycenter} to~\ref{sec:alpha} are strictly conditionally and unconditionally model consistent for the corresponding regressions. For the conditional mean these are the log score $D_{\mathrm{KL}}(y\,\|\,q)$, for categorical and compositional responses alike, and every Bregman loss of the mixture chart, among them the squared loss $\|y - q\|^2$. For the conditional normalized geometric mean they are the Bregman losses of the exponential chart, among them the squared loss $\|\operatorname{alr}(y) - \operatorname{alr}(q)\|^2$ of the log-ratios, for interior responses. For the conditional Fr\'echet mean it is $\tfrac{1}{2}d^2$, for interior and categorical responses. For the conditional power mean $\Pi^{(\alpha)}(F_{Y|Z})$, $\alpha \in (-1, 1)$, it is the restricted $\alpha$-divergence $D^{(\alpha)}[y:q]$.

The plug-in formula of Proposition~\ref{prop:power-mean} does not carry over to regression. Correct specification of $m$ fixes the conditional power mean, and hence the conditional moments only up to a common positive factor, $\E[Y_i^{(1-\alpha)/2} \,|\, Z] = c(Z)\,m_i(Z, \theta_0)^{(1-\alpha)/2}$; the model leaves $c(Z)$ free. The M-estimator maximizes $\sum_k\sum_i y_{k,i}^{(1-\alpha)/2}\,m_i(z_k, \theta)^{(1+\alpha)/2}$. The factor $c(Z)$ only scales the conditional expected loss, so it does not move the minimizer. Regressing the powers $Y_i^{(1-\alpha)/2}$ on the covariates and normalizing estimates the products $c(Z)\,m_i(Z, \theta_0)^{(1-\alpha)/2}$ instead. It is consistent only if these products themselves follow the regression model, an assumption that correct specification of $m$ alone does not imply.

Under correct specification every model estimates the same target, and the choice between models comes down to computation, interpretation and behaviour under misspecification.

A geodesic model of the flat connection whose barycenter is the target is affine in its flat chart $\varphi$. The correctly specified model is then the linear model $\E[\varphi(Y) \,|\, Z] = a + BZ$, and least squares of $\varphi(Y)$ on $Z$ is a consistent M-estimator in closed form. For the exponential connection this is the log-ratio regression of compositional data analysis. Its target is the conditional normalized geometric mean rather than the expected shares, and interior responses, away from the boundary, are required. For the mixture connection it is linear regression of the composition, whose fitted values may lie outside the simplex for covariates outside a bounded range.

Multinomial logistic regression pairs the exponential geodesic model with the mixture barycenter and its canonical divergence, the log score. In the natural parameters, the log score is the convex log-partition function plus a linear term. The empirical risk is therefore convex in $\theta$, the regression form of the uniqueness of m-projections onto exponential families \citep{amari2000}. Its first-order conditions $\sum_k z_k(y_k - m(z_k, \hat\theta)) = 0$ match fitted and observed shares within each covariate, the constant included, and zero shares are allowed. The log-linear form, however, imposes the independence of irrelevant alternatives \citep{mcfadden1974}, and the effects of the coefficients on the shares vary with $z$.

Under misspecification, each consistent loss estimates its own pseudo-true parameter, the minimizer of the expected loss over the model, as \citet{white1982} showed for the likelihood. The loss, which does not affect the target under correct specification, then decides what is estimated.

For a loss in the regularity class of Osband's principle, differentiating $\loss(m(z_k, \theta), y_k)$ with respect to $\theta_j$ by the chain rule and using \eqref{eq:osband-pointwise} yields the estimating equation of the M-estimator,
\begin{equation}\label{eq:model-score}
\sum_k H_{m(z_k, \hat\theta)}(X(m(z_k, \hat\theta), y_k), \partial_jm(z_k, \hat\theta)) = 0, \quad j = 1, \dots, s.
\end{equation}
Here $H_q$ is the Osband tensor of the loss at the fitted composition $q = m(z_k, \hat\theta)$, the bilinear form of Theorem~\ref{thm:osband} that turns the residual into the derivative of the loss, $d_q\loss(\cdot, y) = -H_q(X(q, y), \cdot)$. For the canonical divergence it is the metric $g_q$, and for the Bregman losses of a flat chart it is the Hessian of the generator (Section~\ref{sec:osband}). The vectors $\partial_jm = \partial m/\partial\theta_j$ are the columns of the model differential $d_\theta m$ of Section~\ref{sec:Intro}, the directions in which the fitted composition moves when a parameter changes. The left side of \eqref{eq:model-score} applies the adjoint $(d_\theta m)^{\ast}: T_q^{\ast}\simplex \to \R^s$ of the model differential to the linear forms $H_q(X(q, y_k), \cdot)$ and sums over the observations. These vectors play the role of the covariates in the normal equations: on $\R$, with $m = z\T\theta$, $\partial_jm = z_j$ and $H = 1$, \eqref{eq:model-score} is the normal equations. As on $\R$, these instruments are required for identification. The unconditional expected residual $\E[X(m(Z, \theta), Y)]$ alone does not identify $\theta_0$ in general. \citetalias{dimitriadis2024} (their Example~S1) show this for linear regression with covariates of mean zero, and under a full-rank condition the instruments restore it (their Proposition~S2). The instruments are valid by construction: under correct specification $\E[X(m(Z, \theta_0), Y) \,|\, Z] = 0$, so every function of the covariates is one. With more functions of the covariates than parameters, the test of overidentifying restrictions of \citet{hansen1982} tests correct specification.

In multinomial logistic regression the natural parameters $\eta_i = z\T\theta_i$ are linear in the parameters, so a change $d\eta$ moves the fitted composition by $dq_i = q_i(d\eta_i - \sum_j q_j\,d\eta_j)$. With the log score, whose Osband tensor is the Fisher-Rao metric, the pairing is $g_q(y - q, dq) = \sum_i (y_i - q_i)\,dq_i/q_i = \sum_i (y_i - q_i)\,d\eta_i$, because $\sum_i (y_i - q_i) = 0$. The Fisher weights $1/q_i$ cancel against the derivative of the exponential chart, and \eqref{eq:model-score} becomes $\sum_k z_k(y_{k,i} - m_i(z_k, \hat\theta)) = 0$ for every category $i$, with the covariates themselves as instruments. This is the canonical link of generalized linear models, seen as the duality of the exponential and mixture connections.

A gauge transformation $A$ of the residual, with $A_q$ invertible, leaves the barycenter unchanged, since $A_q\bar X(q, F) = 0$ exactly when $\bar X(q, F) = 0$. It does change the estimating equation of a model,
\[
\Psi_A(\theta)_j = g_{q(\theta)}(A_{q(\theta)}\bar X(q(\theta), F_z), \partial_jm(z, \theta)) = 0, \quad q(\theta) = m(z, \theta), \quad F_z = F_{Y|Z = z},
\]
since $A_q$ is in general not an isometry, so different gauges produce different method-of-moments estimators of $\theta$. Under correct specification all of them are solved at $\theta_0$, where $\bar X$ vanishes. Differentiating $\Psi_A$ at $\theta_0$, every term except the one that differentiates $\bar X$ contains the factor $\bar X(q(\theta_0), F_z) = 0$, so the Jacobian is $G_A(u, v) = g(A\nabla_{d_\theta m\,u}\bar X, d_\theta m\,v)$. For the gauge $h$ of a consistent loss, with $H = g(h\,\cdot, \cdot)$, \eqref{eq:osband-simplex} turns $G_h$ into minus the Hessian in $\theta$ of the conditional expected loss. $G_h$ is symmetric, and $\Psi_h$ is minus the gradient of that loss. For other gauges $G_A$ is in general not symmetric, so $\Psi_A$ is not the gradient of any function of $\theta$. Its equation is then not the first-order condition of a minimization, as for the barycenter equation in case (iii) of Corollary~\ref{cor:regimes}. At $\alpha = 0$ the multiples of the identity are the only gauges of consistent losses, by the same corollary.

The influence function of an M-estimator, viewed as the functional $T(F) = \arg\min_q \bar\loss(q, F)$ of the distribution \citep{huber1964}, is the derivative of $T$ at $F$ towards a point mass at $x$, $\mathrm{IF}(x) = \lim_{\varepsilon \downarrow 0}\varepsilon^{-1}\left(T((1 - \varepsilon)F + \varepsilon\delta_x) - T(F)\right)$, and its supremum over $x$ is the gross-error sensitivity. On a convex class of distributions all strictly consistent losses for a functional define the same $T$ and hence the same influence function; for the mean it is $x - \E_F X$, unbounded in the observation. In a regression model the contaminated distributions lie outside the model, each loss has its own pseudo-true parameter there, and the influence function depends on the loss. Differentiating the estimating equation along the contamination yields the inverse Hessian of the expected loss applied to the summand of \eqref{eq:model-score} at the contaminating observation $(z, y)$, so the influence of an observation is bounded in the response exactly when that summand is; for linear regression the summand $z(y - z\T\theta)$ is unbounded in $y$.

The summand $H_q(X(q, y), \partial_jm)$, with $H_q$ continuous in $q$, is bounded in $y$ whenever the residual is. This is the case for the mixture residual $y - q$, for the Levi-Civita residual $\log_q(y)$, of length $d(q, y) < \pi$, and for the score of the power-mean losses, the projected cone residual $P_qX^{(\alpha)}(q, y)$, whose chart $m \mapsto m^{(1 - \alpha)/2}$ maps the simplex into the unit cube; it is not the case for the exponential residual, which contains $\log(y_i/q_i)$. Multinomial logistic regression, whose summand $z(y_i - m_i(z, \theta))$ is bounded in the response, is the mixture case. Without covariates, the influence function of the Fr\'echet mean is $\Hess\bar\loss(q^\ast, F)^{-1}\log_{q^\ast}(x)$ at $q^\ast = T(F)$, of length at most $\pi/c_F$, where $c_F = \E_F[(d(q^\ast, X)/2)\cot(d(q^\ast, X)/2)]$ averages the transverse eigenvalue of the Hessian of Section~\ref{sec:barycenter} and is positive for every $F \in \Ff_0$, since $d(q^\ast, X) < \pi$ \citep[Section~3.6.4]{navasconvex}. The bound comes from the bounded residual, and curvature enters only through $c_F$: the mixture connection is flat with a bounded residual, the exponential connection flat with an unbounded one. Whether $c_F$ is bounded below over $\Ff_0$ is open.

The estimating equation \eqref{eq:model-score} is a method-of-moments equation with the instruments $\partial_jm$, as in the generalized method of moments of \citet{hansen1982} with as many equations as parameters, and the theory of such estimators yields the limit distribution of $\hat\theta$. Write $\psi_\theta(z, y) \in \R^s$ for its summand, $\psi_\theta(z, y)_j = H_{m(z, \theta)}(X(m(z, \theta), y), \partial_jm(z, \theta))$, or for the summand of $\Psi_A$ with another gauge $A$. Under correct specification $\E\psi_{\theta_0}(Z, Y) = 0$, and the Jacobian of $\theta \mapsto \E\psi_\theta(Z, Y)$ at $\theta_0$ is the expectation $G$ over $Z$ of the matrix $G_A$ above. If $G$ is nonsingular, which is the rank condition expressing the relevance of the instruments, and $\psi_\theta$ is Lipschitz in $\theta$ near $\theta_0$ with a square-integrable constant, a consistent root $\hat\theta$ satisfies $\sqrt{N}(\hat\theta - \theta_0) \to N(0, G^{-1}\E[\psi_{\theta_0}\psi_{\theta_0}\T]G^{-\top})$ \citep[Theorem~5.21]{vandervaart1998}. For the gauge of a consistent loss $G$ is symmetric, and the limit is the sandwich covariance of the M-estimator. For the other gauges, and for the barycenter equation of case (iii) of Corollary~\ref{cor:regimes}, the limit has the same form although no loss is minimized: the method of moments requires identification only. Where the residual is bounded, as for the mixture, Levi-Civita and power-mean losses, the second moment of $\psi$ is finite as soon as that of $d_\theta m$ is, and no moment of the response is required. Without covariates every consistent loss and every gauge yield the same estimator, the root of $\sum_k X(\hat q, x_k) = 0$, and the same limit. On a class $\Ff_B$ this root is unique, since the empirical distribution lies in $\Ff_B$, and it is consistent when $G$ is nonsingular, because a nonsingular $G$ places a root near the barycenter with probability tending to one \citep[Theorem~5.42]{vandervaart1998}. For the Fr\'echet mean $G$ is minus the Hessian of the expected loss, nonsingular since the Hessian is at least $c_F\,g$, and the score is bounded by $\pi$. Consistency follows from Theorem~\ref{thm:lcf-consistent} and the continuity of the loss on the closed simplex, which is compact \citep[Theorem~5.7]{vandervaart1998}. So, for every $F \in \Ff_0$ and without moment conditions, $\sqrt{N}\log_{q^\ast}(\hat q)$ is asymptotically normal with the covariance $\Hess\bar\loss(q^\ast, F)^{-1}\,\E_F[\log_{q^\ast}(X) \otimes \log_{q^\ast}(X)]\,\Hess\bar\loss(q^\ast, F)^{-1}$ of the influence function. In particular the sample Fr\'echet mean on the simplex is not smeary in the sense of \citet{eltznerhuckemann2019}, whose slower rates arise when the Hessian degenerates.

\section{Calibration and discrimination of compositional forecasts}\label{sec:calibration}

On the real line, a forecast $Q$ of an outcome $Y$ is scored by the squared loss $(Q - Y)^2$, which is consistent for the mean. The forecast is \emph{auto-calibrated} if $\E[Y \,|\, Q] = Q$, that is, if the outcomes are on average what the forecast announces. With the recalibrated forecast $\tilde Q = \E[Y \,|\, Q]$, which corrects the forecast by conditioning on it, and the unconditional average $\bar y = \E Y$, the expected score splits as $\E(Q - Y)^2 = \E(Y - \bar y)^2 - \E(\tilde Q - \bar y)^2 + \E(Q - \tilde Q)^2$: the uncertainty of the outcome, minus the discrimination, the variation of the recalibrated forecast, plus the miscalibration, the expected squared distance between the issued and the recalibrated forecast.

On the simplex the mean is replaced by a functional $\Gamma$, a barycenter or a power mean, and the squared loss by a loss strictly consistent for it. Let $Q$ be the issued forecast, a random point of the simplex, $\tilde Q = \Gamma(F_{Y|Q})$ the recalibrated forecast and $\bar q = \Gamma(F_Y)$ the unconditional average. The forecast is auto-calibrated if $\tilde Q = Q$ almost surely, that is, if $\E[V(Q, Y) \,|\, Q] = 0$ for an identification function $V$ of $\Gamma$, a conditional moment condition. For the Fr\'echet mean this says that, given the forecast, the Fr\'echet mean of the outcome is the forecast itself: the outcomes are centered on the forecast in the Fisher-Rao geometry, as they are centered on it in the Euclidean sense on the real line. For every loss $\loss$ strictly consistent for $\Gamma$ on a class containing the conditional laws $F_{Y|Q}$, the expected score splits in the same way,
\begin{equation}\label{eq:score-decomposition}
\E\,\loss(Q, Y) = \mathrm{UNC} - \mathrm{DSC} + \mathrm{MCB},
\end{equation}
with $\mathrm{UNC} = \E\,\loss(\bar q, Y)$, $\mathrm{DSC} = \E\,\loss(\bar q, Y) - \E\,\loss(\tilde Q, Y)$ and $\mathrm{MCB} = \E\,\loss(Q, Y) - \E\,\loss(\tilde Q, Y)$. On the real line with the squared loss, these are the three terms above. Conditionally on $Q$ the recalibrated forecast minimizes the expected loss, so $\mathrm{DSC} \ge 0$ and $\mathrm{MCB} \ge 0$, with $\mathrm{MCB} = 0$ exactly for auto-calibrated forecasts. The argument uses only strict consistency, and the decomposition that \citet*{dimitriadisgneitingjordan2021} use for probabilistic classifiers therefore holds for every functional and loss of this paper. Calibration matters wherever compositions are produced by classifiers: \citet{fernandesmachado2025}, for instance, obtain compositions from multinomial logistic regression and name a poorly calibrated classifier as a source of bias in their counterfactuals.

On the real line, miscalibration and discrimination are expected squared distances between forecasts. This holds again at the ends $\alpha = \pm 1$ of the $\alpha$-family, where the connections of the simplex are flat: the mixture connection at $\alpha = -1$, whose barycenter is the arithmetic mean, and the exponential connection at $\alpha = 1$, whose barycenter is the normalized geometric mean. For a flat connection with flat chart $\varphi$ and a Bregman loss $\loss(q, y) = D_\psi[\varphi(y):\varphi(q)]$ the barycenter is the mean in the chart (Section~\ref{sec:osband}), and the two-stage identity \eqref{eq:two-stage}, applied in the chart conditionally on $Q$ with the arguments $Q$ and $\bar q$, yields
\[
\mathrm{MCB} = \E\,D_\psi[\varphi(\tilde Q):\varphi(Q)], \qquad \mathrm{DSC} = \E\,D_\psi[\varphi(\tilde Q):\varphi(\bar q)].
\]
Miscalibration then depends on the pair of issued and recalibrated forecasts alone, and not on how the outcomes spread around the recalibrated forecast. The reliability diagram, the plot of $\tilde Q$ against $Q$, therefore determines the miscalibration \citep{dimitriadisgneitingjordan2021}. At $\alpha = -1$, with the log score, $\mathrm{MCB} = \E\,D_{\mathrm{KL}}(\tilde Q\,\|\,Q)$ and $\mathrm{DSC} = \E\,D_{\mathrm{KL}}(\tilde Q\,\|\,\bar q)$; at $\alpha = 1$, with the loss $D_{\mathrm{KL}}(q\,\|\,y)$ and interior outcomes, $\mathrm{MCB} = \E\,D_{\mathrm{KL}}(Q\,\|\,\tilde Q)$.

For $\alpha \in (-1, 1)$, with target the power mean $\Pi^{(\alpha)}$ and loss the restricted $\alpha$-divergence $D^{(\alpha)}[y:q]$, write $m_i = \E[Y_i^{(1-\alpha)/2} \,|\, Q]$ for the conditional moments and
\[
c_\alpha(F_{Y|Q}) = \Bigl(\sum_j m_j^{2/(1-\alpha)}\Bigr)^{(1-\alpha)/2}.
\]

\begin{proposition}\label{prop:calibration-alpha}
For $\alpha \in (-1, 1)$ let $c_\alpha = c_\alpha(F_{Y|Q})$. Then $\mathrm{MCB} = \E[c_\alpha D^{(\alpha)}[\tilde Q:Q]]$ and $\mathrm{DSC} = \E[c_\alpha D^{(\alpha)}[\tilde Q:\bar q]]$, where $c_\alpha \in (0, 1]$, with $c_\alpha = 1$ exactly when $Y$ is determined by $Q$.
\end{proposition}

\begin{proof}
Conditionally on $Q$, the expected loss at $q$ is $\frac{4}{1-\alpha^2}(1 - \sum_i m_iq_i^{(1+\alpha)/2})$. By Proposition~\ref{prop:power-mean} its minimizer $\tilde Q$ satisfies $m_i = c_\alpha\,\tilde Q_i^{(1-\alpha)/2}$, so the excess over the minimum at $q$ is $c_\alpha D^{(\alpha)}[\tilde Q:q]$. Taking $q = Q$ and $q = \bar q$ and integrating yields the identities. Jensen's inequality for the concave power $(1-\alpha)/2$ yields $m_j^{2/(1-\alpha)} \le \E[Y_j \,|\, Q]$, with equality exactly when $Y_j$ is determined by $Q$, and summing over $j$ yields $c_\alpha \le 1$.
\end{proof}

Miscalibration is thus still a divergence between the issued and the recalibrated forecast, multiplied by a factor that measures how concentrated the outcome is given the forecast. The concentration is controlled by the conditional dispersion of the outcome in the Fisher-Rao metric: when it is small, $1 - c_\alpha(F_{Y|Q}) \approx \frac{1-\alpha^2}{8}\sum_j \mathrm{Var}(Y_j \,|\, Q)/\E[Y_j \,|\, Q]$, and the sum on the right is $\E[\|Y - \E[Y \,|\, Q]\|_g^2 \,|\, Q]$, the expected squared length of the deviation of the outcome from its conditional mean, measured by the metric \eqref{eq:fr-metric} at that mean. Two forecasters with the same reliability diagram but different conditional dispersion of the outcome therefore incur different miscalibration. The reliability diagram determines the miscalibration at the ends $\alpha = \pm 1$ of the family, where the factor $1 - \alpha^2$ vanishes and $c_\alpha \to 1$, but not in between.

For the Fr\'echet mean with the loss $\tfrac{1}{2}d^2$, the curvature of Section~\ref{sec:residuals} replaces the factor by two bounds. Let $\rho < \pi$ be the diameter of a geodesically convex region containing the outcomes, the forecasts and the recalibrated forecasts, and $c_\rho = (\rho/2)\cot(\rho/2) \in (0, 1)$.

\begin{proposition}\label{prop:calibration-frechet}
$c_\rho\,\E\,\tfrac{1}{2}d(\tilde Q, Q)^2 \le \mathrm{MCB} \le \E\,\tfrac{1}{2}d(\tilde Q, Q)^2$, and the same bounds hold for $\mathrm{DSC}$ with $d(\tilde Q, \bar q)$ in place of $d(\tilde Q, Q)$.
\end{proposition}

\begin{proof}
Conditionally on $Q$, let $R(q) = \E[\tfrac{1}{2}d(q, Y)^2 \,|\, Q]$ and $\gamma$ the geodesic from $\tilde Q$ to $Q$, of speed $d(\tilde Q, Q)$. Since $\grad R(\tilde Q) = 0$, Taylor's formula along $\gamma$ yields $R(Q) - R(\tilde Q) = \int_0^1 (1-t)\Hess R(\gamma(t))(\dot\gamma(t), \dot\gamma(t))\,dt$. The eigenvalues of the Hessian of $\tfrac{1}{2}d(\cdot, y)^2$ at distance $\delta$ from $y$ are $1$ and $(\delta/2)\cot(\delta/2)$ (Section~\ref{sec:barycenter}), both in $[c_\rho, 1]$ for $\delta \le \rho$, since $u\cot u$ decreases on $(0, \pi/2)$. Averaging over $Y$ and integrating in $t$ yields the bounds.
\end{proof}

On the real line, the curvature is zero and $c_\rho = 1$, so miscalibration equals half the expected squared distance; on the simplex it is determined by that distance up to the factor $c_\rho$, which depends only on the curvature and the diameter. Near the vertices, where $\rho$ approaches $\pi$, the factor $c_\rho$ tends to zero, and the lower bound degenerates.

\section{Ranking compositional forecasts}\label{sec:ranking}

Two forecasters issue, at each time $t$, point forecasts $q_t$ and $q'_t$ of the same functional of the conditional law of an outcome $Y_t$, before $Y_t$ is observed. With a strictly consistent loss $\loss$ they are compared through the score differences $\Delta_t = \loss(q_t, Y_t) - \loss(q'_t, Y_t)$, and the first forecaster is better when the differences are negative on average. Strict consistency means that, for every conditional law, the expected score is minimized exactly at the value of the functional, and many losses have this property; for the mean on the real line, every Bregman loss does. \citet{patton2020} shows that two of them can rank the same pair of forecasters differently as soon as a forecast comes from a misspecified model, contains estimation error, or uses information not nested in that of the other forecaster. Only correctly specified forecasts built on nested information are ranked in the same way by all consistent losses. In general, which forecaster is better depends on the loss chosen for the comparison.

On the simplex the functional itself is a matter of choice. The functional to be forecast is the barycenter of a connection, and choosing the connection chooses the functional: the arithmetic mean for the mixture connection, the normalized geometric mean for the exponential connection, the Fr\'echet mean for the Levi-Civita connection of the Fisher-Rao metric. The choice is problem-specific, and there is no preferred data-driven way to make it. The geometry determines the functional, and Corollary~\ref{cor:regimes} then determines its consistent losses in the regularity class of Osband's principle. At the flat ends $\alpha = \pm 1$ the consistent losses are the Bregman losses of all strictly convex generators in the flat chart, and the ranking of forecasters depends on the generator, as on the real line. At $\alpha = 0$, with $n \ge 3$ categories, the consistent losses are the multiples of $\tfrac{1}{2}d^2$, up to a function of the outcome, and all of them rank the forecasters in the same way. For $\alpha \notin \{-1, 0, 1\}$ no consistent loss exists, and forecasts of these barycenters cannot be ranked by any score. The residual still identifies the barycenter, so whether a forecaster issues it can be tested through the conditional moment condition $\E[X(Q, Y) \,|\, Q] = 0$. In the Fisher-Rao geometry the choice of the functional is therefore the only choice: once it is made, consistency alone decides which forecaster is better.

For the mixture connection, the case $\alpha = -1$, the functional is the arithmetic mean, and two of its consistent losses are used frequently. One is the quadratic loss $\|y - q\|^2$, which is the Brier score when the outcomes are vertices. The other is the Kullback-Leibler loss $D_{\mathrm{KL}}(y\,\|\,q)$, which is the logarithmic score up to a function of the outcome. They are the Bregman losses of the generators $\|p\|^2$ and $\sum_i p_i\log p_i$. For a forecast $q$ that does not change over time, the two-stage identity \eqref{eq:two-stage} yields $\E D_\psi[Y:q] = \E D_\psi[Y:\mu] + D_\psi[\mu:q]$ with $\mu = \E Y$. Two such forecasts are therefore ranked by the excess $D_\psi[\mu:q]$ alone.

\begin{table}[t]
\caption{Excess expected scores of two constant forecasts of the arithmetic mean $\mu = (0.45, 0.45, 0.10)$.}\label{tab:reversal}
\begin{tabular}{llcc}
\hline
forecast & components & $\|\mu - q\|^2$ & $D_{\mathrm{KL}}(\mu\,\|\,q)$ \\
\hline
$q$ & $(0.40, 0.40, 0.20)$ & $0.0150$ & $0.0367$ \\
$q'$ & $(0.49, 0.49, 0.02)$ & $0.0096$ & $0.0843$ \\
\hline
\end{tabular}
\end{table}

\begin{example}
Let the outcome have mean $\mu = (0.45, 0.45, 0.10)$, and let two forecasters issue the constant forecasts $q$ and $q'$ of Table~\ref{tab:reversal}, neither of which equals $\mu$. The quadratic loss prefers $q'$, and the Kullback-Leibler loss prefers $q$. The quadratic loss measures errors in the coordinates and rewards $q'$ for being close on the two large categories. To second order, the Kullback-Leibler loss weighs the squared error in category $i$ by $1/q_i$, the weight of the Fisher-Rao metric \eqref{eq:fr-metric}, and punishes $q'$ for assigning $0.02$ to a category whose mean share is $0.10$.
\end{example}

At $\alpha = 0$ the functional is the Fr\'echet mean. By Corollary~\ref{cor:regimes}~(ii), every consistent loss in the regularity class of Osband's principle has the form $\loss(q, y) = \tfrac{c}{2}d(q, y)^2 + a(y)$ with $c > 0$. The function $a$ of the outcome cancels in the score difference, $\Delta_t = \tfrac{c}{2}(d(q_t, Y_t)^2 - d(q'_t, Y_t)^2)$, and the constant $c$ only rescales it. Every consistent loss therefore yields the same sign of $\Delta_t$ at every time, and hence the same ranking of the two forecasters, whatever the forecasts and whatever the law of the outcomes. Unlike on the real line, this requires neither correct specification nor nested information. In the terminology of \citet*{ehm2016}, a forecaster who is better under one consistent loss dominates, that is, is better under all of them, and no Murphy diagram is required to check it. The condition $n \ge 3$ is necessary, as in Theorem~\ref{thm:rigidity}. For $n = 2$ the simplex is the segment of binary probabilities, and the chart $p \mapsto 2\arcsin\sqrt{p_1}$ maps it isometrically onto an interval. The Fr\'echet mean is the mean in this chart, and every Bregman loss of the chart is consistent, so rankings can again depend on the loss: with two categories there is no curvature, and one recovers the case of the real line.

\section*{Acknowledgments}
The author used Claude Opus 5, a large language model by Anthropic, to write the script of the numerical check in the Supplementary Material, to shorten the proofs and to check the references. The mathematical content is the author's, and all AI-assisted material was verified by the author.

\bibliographystyle{plainnat}
\bibliography{references}

\begin{thebibliography}{31}
\providecommand{\natexlab}[1]{#1}
\providecommand{\url}[1]{\texttt{#1}}
\expandafter\ifx\csname urlstyle\endcsname\relax
  \providecommand{\doi}[1]{doi: #1}\else
  \providecommand{\doi}{doi: \begingroup \urlstyle{rm}\Url}\fi

\bibitem[Amari and Nagaoka(2000)]{amari2000}
S.~Amari and H.~Nagaoka.
\newblock \emph{Methods of Information Geometry}, volume 191 of
  \emph{Translations of Mathematical Monographs}.
\newblock American Mathematical Society / Oxford University Press, 2000.

\bibitem[Ay and Amari(2015)]{ayamari2015}
N.~Ay and S.~Amari.
\newblock A novel approach to canonical divergences within information
  geometry.
\newblock \emph{Entropy}, 17\penalty0 (12):\penalty0 8111--8129, 2015.

\bibitem[Ay et~al.(2017)Ay, Jost, L{\^e}, and Schwachh{\"o}fer]{ayjost2017}
N.~Ay, J.~Jost, H.~V. L{\^e}, and L.~Schwachh{\"o}fer.
\newblock \emph{Information Geometry}.
\newblock Ergebnisse der Mathematik und ihrer Grenzgebiete. Springer, Cham,
  2017.

\bibitem[Banerjee et~al.(2005)Banerjee, Merugu, Dhillon, and
  Ghosh]{banerjee2005}
A.~Banerjee, S.~Merugu, I.~S. Dhillon, and J.~Ghosh.
\newblock Clustering with {B}regman divergences.
\newblock \emph{Journal of Machine Learning Research}, 6:\penalty0 1705--1749,
  2005.

\bibitem[Bhattacharya and Patrangenaru(2003)]{bhattacharya2003}
R.~Bhattacharya and V.~Patrangenaru.
\newblock Large sample theory of intrinsic and extrinsic sample means on
  manifolds {I}.
\newblock \emph{The Annals of Statistics}, 31:\penalty0 1--29, 2003.

\bibitem[Dimitriadis et~al.(2021)Dimitriadis, Gneiting, and
  Jordan]{dimitriadisgneitingjordan2021}
T.~Dimitriadis, T.~Gneiting, and A.~I. Jordan.
\newblock Stable reliability diagrams for probabilistic classifiers.
\newblock \emph{Proceedings of the National Academy of Sciences}, 118:\penalty0
  e2016191118, 2021.

\bibitem[Dimitriadis et~al.(2024{\natexlab{a}})Dimitriadis, Fissler, and
  Ziegel]{dfzidentification}
T.~Dimitriadis, T.~Fissler, and J.~F. Ziegel.
\newblock Osband's principle for identification functions.
\newblock \emph{Statistical Papers}, 65:\penalty0 1125--1132,
  2024{\natexlab{a}}.

\bibitem[Dimitriadis et~al.(2024{\natexlab{b}})Dimitriadis, Fissler, and
  Ziegel]{dimitriadis2024}
T.~Dimitriadis, T.~Fissler, and J.~F. Ziegel.
\newblock Characterizing {M}-estimators.
\newblock \emph{Biometrika}, 111:\penalty0 339--346, 2024{\natexlab{b}}.

\bibitem[Eguchi(1992)]{eguchi1992}
S.~Eguchi.
\newblock Geometry of minimum contrast.
\newblock \emph{Hiroshima Mathematical Journal}, 22\penalty0 (3):\penalty0
  631--647, 1992.

\bibitem[Eguchi and Komori(2022)]{eguchikomori2022}
S.~Eguchi and O.~Komori.
\newblock \emph{Minimum Divergence Methods in Statistical Machine Learning:
  From an Information Geometric Viewpoint}.
\newblock Springer, Tokyo, 2022.

\bibitem[Ehm et~al.(2016)Ehm, Gneiting, Jordan, and Kr{\"u}ger]{ehm2016}
W.~Ehm, T.~Gneiting, A.~Jordan, and F.~Kr{\"u}ger.
\newblock Of quantiles and expectiles: consistent scoring functions, {C}hoquet
  representations and forecast rankings.
\newblock \emph{Journal of the Royal Statistical Society: Series B},
  78:\penalty0 505--562, 2016.

\bibitem[Eltzner and Huckemann(2019)]{eltznerhuckemann2019}
B.~Eltzner and S.~F. Huckemann.
\newblock A smeary central limit theorem for manifolds with application to
  high-dimensional spheres.
\newblock \emph{The Annals of Statistics}, 47:\penalty0 3360--3381, 2019.

\bibitem[Fernandes~Machado et~al.(2025)Fernandes~Machado, Gallic, and
  Charpentier]{fernandesmachado2025}
A.~Fernandes~Machado, E.~Gallic, and A.~Charpentier.
\newblock Optimal transport on categorical data for counterfactuals using
  compositional data and {D}irichlet transport.
\newblock In J.~Kwok, editor, \emph{Proceedings of the Thirty-Fourth
  International Joint Conference on Artificial Intelligence}, pages 5136--5144.
  International Joint Conferences on Artificial Intelligence Organization,
  2025.
\newblock \doi{10.24963/ijcai.2025/572}.

\bibitem[Fissler and Ziegel(2016)]{fisslerziegel2016}
T.~Fissler and J.~F. Ziegel.
\newblock Higher order elicitability and {O}sband's principle.
\newblock \emph{The Annals of Statistics}, 44:\penalty0 1680--1707, 2016.

\bibitem[Fletcher(2013)]{fletcher2013}
P.~T. Fletcher.
\newblock Geodesic regression and the theory of least squares on {R}iemannian
  manifolds.
\newblock \emph{International Journal of Computer Vision}, 105:\penalty0
  171--185, 2013.

\bibitem[Gneiting(2011)]{gneiting2011}
T.~Gneiting.
\newblock Making and evaluating point forecasts.
\newblock \emph{Journal of the American Statistical Association}, 106:\penalty0
  746--762, 2011.

\bibitem[Hansen(1982)]{hansen1982}
L.~P. Hansen.
\newblock Large sample properties of generalized method of moments estimators.
\newblock \emph{Econometrica}, 50:\penalty0 1029--1054, 1982.

\bibitem[Huber(1964)]{huber1964}
P.~J. Huber.
\newblock Robust estimation of a location parameter.
\newblock \emph{Annals of Mathematical Statistics}, 35\penalty0 (1):\penalty0
  73--101, 1964.

\bibitem[Karcher(1977)]{karcher1977}
H.~Karcher.
\newblock {Riemannian} center of mass and mollifier smoothing.
\newblock \emph{Communications on Pure and Applied Mathematics}, 30\penalty0
  (5):\penalty0 509--541, 1977.

\bibitem[Kobayashi and Ohno(2022)]{kobayashiohno2022}
S.~Kobayashi and Y.~Ohno.
\newblock On a constant curvature statistical manifold.
\newblock \emph{Information Geometry}, 5:\penalty0 31--46, 2022.

\bibitem[Kurose(1990)]{kurose1990}
T.~Kurose.
\newblock Dual connections and affine geometry.
\newblock \emph{Mathematische Zeitschrift}, 203:\penalty0 115--121, 1990.

\bibitem[McFadden(1974)]{mcfadden1974}
D.~McFadden.
\newblock Conditional logit analysis of qualitative choice behavior.
\newblock In P.~Zarembka, editor, \emph{Frontiers in Econometrics}, pages
  105--142. Academic Press, New York, 1974.

\bibitem[Navas(2026{\natexlab{a}})]{navasconvex}
P.~Navas.
\newblock Estimators from convex geometries of the simplex.
\newblock arXiv:2610.21551, 2026{\natexlab{a}}.

\bibitem[Navas(2026{\natexlab{b}})]{navassupp}
P.~Navas.
\newblock Supplement to ``{C}haracterizing {M}-estimators using canonical
  divergences from information geometry''.
\newblock Appended to this version, 2026{\natexlab{b}}.

\bibitem[Osband(1985)]{osband1985}
K.~H. Osband.
\newblock \emph{Providing Incentives for Better Cost Forecasting}.
\newblock PhD thesis, University of California, Berkeley, 1985.

\bibitem[Patton(2020)]{patton2020}
A.~J. Patton.
\newblock Comparing possibly misspecified forecasts.
\newblock \emph{Journal of Business \& Economic Statistics}, 38:\penalty0
  796--809, 2020.

\bibitem[Pennec(2019)]{pennec2019}
X.~Pennec.
\newblock Curvature effects on the empirical mean in {R}iemannian and affine
  manifolds: a non-asymptotic high concentration expansion in the small-sample
  regime.
\newblock arXiv:1906.07418, 2019.

\bibitem[Pennec and Arsigny(2013)]{pennec2013}
X.~Pennec and V.~Arsigny.
\newblock Exponential barycenters of the canonical {C}artan connection and
  invariant means on {L}ie groups.
\newblock In F.~Nielsen and R.~Bhatia, editors, \emph{Matrix Information
  Geometry}, pages 123--166. Springer, Berlin, Heidelberg, 2013.

\bibitem[Savage(1971)]{savage1971}
L.~J. Savage.
\newblock Elicitation of personal probabilities and expectations.
\newblock \emph{Journal of the American Statistical Association}, 66:\penalty0
  783--801, 1971.

\bibitem[van~der Vaart(1998)]{vandervaart1998}
A.~W. van~der Vaart.
\newblock \emph{Asymptotic Statistics}.
\newblock Cambridge University Press, Cambridge, 1998.

\bibitem[White(1982)]{white1982}
H.~White.
\newblock Maximum likelihood estimation of misspecified models.
\newblock \emph{Econometrica}, 50:\penalty0 1--25, 1982.

\end{thebibliography}

\clearpage
\setcounter{section}{0}
\renewcommand{\thesection}{S\arabic{section}}
\renewcommand{\theHsection}{S\arabic{section}}
\section*{Supplementary Material}
This supplement contains the proofs of Lemma~\ref{lem:affine-factorization} and Theorem~\ref{thm:lcf-consistent} and a numerical check of Proposition~\ref{prop:nonintegrable}. Numbers of the form 4.4 refer to the paper, numbers of the form S1.1 to this supplement.

\section{Affine factorization: proof of Lemma 4.4}\label{sec:S-lemma}

Lemma~4.4 of the paper is used for Proposition~4.5, which shows that every strict identification function of a barycenter is a gauge transformation of the residual, and for Theorem~5.1, Osband's principle on the simplex. We restate it.

\begin{lemma}[Lemma 4.4]\label{lem:S-affine}
Let $\Ff$ be convex, let $V$ and $W$ be finite-dimensional vector spaces, and let $\Phi: \Ff \to V$ and $\Psi: \Ff \to W$ be affine maps such that $0$ is an interior point of the convex set $\Phi(\Ff)$ and $\Psi(F) = 0$ whenever $\Phi(F) = 0$. Then there is a linear map $L: V \to W$ with $\Psi = L \circ \Phi$. If moreover $\Psi(F) = 0$ only when $\Phi(F) = 0$, then $L$ is injective.
\end{lemma}

\begin{proof}
Affine means that $\Phi((1 - \lambda)F_0 + \lambda F_1) = (1 - \lambda)\Phi(F_0) + \lambda\Phi(F_1)$ for $F_0, F_1 \in \Ff$ and $\lambda \in [0, 1]$, and likewise for $\Psi$. Write $C = \Phi(\Ff)$. It is convex, since $\Ff$ is convex and $\Phi$ affine, and by hypothesis it contains a ball $\{v \in V : |v| \le \varepsilon\}$ for some norm $|\cdot|$ on $V$ and some $\varepsilon > 0$.

The first step shows that $\Psi(F)$ depends on $F$ only through $\Phi(F)$. Let $F, F' \in \Ff$ with $\Phi(F) = \Phi(F') = v$. If $v = 0$, both $\Psi(F)$ and $\Psi(F')$ vanish by hypothesis. If $v \ne 0$, choose $s > 0$ with $-sv \in C$, which exists since $0$ is interior to $C$, and $F_1 \in \Ff$ with $\Phi(F_1) = -sv$. With $\lambda = s/(1 + s)$ the distribution $G = \lambda F + (1 - \lambda)F_1$ lies in $\Ff$, and
\[
\Phi(G) = \lambda v - (1 - \lambda)sv = \frac{s}{1 + s}\,v - \frac{s}{1 + s}\,v = 0,
\]
so $\Psi(G) = 0$, that is, $\lambda\Psi(F) + (1 - \lambda)\Psi(F_1) = 0$. The same identity holds with $F'$ in place of $F$, and subtracting the two yields $\lambda(\Psi(F) - \Psi(F')) = 0$, hence $\Psi(F) = \Psi(F')$ since $\lambda > 0$. So $L_0(v) = \Psi(F)$, for any $F \in \Ff$ with $\Phi(F) = v$, defines a map $L_0: C \to W$ with $\Psi = L_0 \circ \Phi$ and $L_0(0) = 0$.

The second step shows that $L_0$ is affine on $C$. For $v, v' \in C$ choose $F, F' \in \Ff$ with $\Phi(F) = v$ and $\Phi(F') = v'$. For $\lambda \in [0, 1]$,
\[
L_0((1 - \lambda)v + \lambda v') = \Psi((1 - \lambda)F + \lambda F') = (1 - \lambda)L_0(v) + \lambda L_0(v').
\]
With $v' = 0$ this yields $L_0(\mu v) = \mu L_0(v)$ for $v \in C$ and $\mu \in [0, 1]$. With $\lambda = 1/2$, and whenever $v$, $v'$ and $v + v'$ lie in $C$, it yields $L_0(v + v') = 2L_0(\tfrac{1}{2}(v + v')) = L_0(v) + L_0(v')$.

The third step extends $L_0$ to $V$. For $v \in V$ choose $\mu > 0$ with $\mu v \in C$, which exists since $C$ contains a ball around $0$, and set $L(v) = L_0(\mu v)/\mu$. This does not depend on $\mu$: for $0 < \mu < \mu'$ with $\mu' v \in C$, homogeneity yields $L_0(\mu v) = L_0((\mu/\mu')\,\mu' v) = (\mu/\mu')L_0(\mu' v)$. For $v \in C$ one may take $\mu = 1$, so $L$ extends $L_0$ and $\Psi = L \circ \Phi$. Positive homogeneity of $L$ follows from its definition, additivity from that of $L_0$ after scaling $v$, $v'$ and $v + v'$ into the ball, and $L(-v) = -L(v)$ from $L_0(\mu v) + L_0(-\mu v) = L_0(0) = 0$ for small $\mu > 0$. So $L$ is linear.

For the last claim assume that $\Psi(F) = 0$ only when $\Phi(F) = 0$, and let $L(v) = 0$. Choose $\mu > 0$ with $\mu v \in C$ and $F \in \Ff$ with $\Phi(F) = \mu v$. Then $\Psi(F) = L(\mu v) = \mu L(v) = 0$, so $\Phi(F) = 0$, that is, $\mu v = 0$ and $v = 0$.
\end{proof}

\begin{remark}\label{rem:S-generality}
The lemma covers more than its two uses in the paper. The map $L$ is unique, since $C$ contains a ball and therefore spans $V$, and two linear maps that agree on $C$ agree on $V$. The proof does not use that $W$ has finite dimension, so $\Psi$ may take values in any real vector space. Of $V$ it uses only that every $v \in V$ has a multiple $\mu v$ in $C$, which holds when $V$ is a topological vector space and $0$ is interior to $C$. Of the class $\Ff$ only convexity is used, and of $\Phi$ only that $0$ is interior to $\Phi(\Ff)$, which for $\Phi = \bar X(q, \cdot)$ is Assumption~4.3 (V1). Proposition~4.5 and Theorem~5.1 therefore hold on every convex class on which (V1) holds, not only on $\Ff_0$ and $\Ff_B$. Classes of laws with a density, or with a bound on a moment, are convex, and (V1) holds on them whenever they contain, for each $q$ in the image of the barycenter, laws whose expected residuals at $q$ surround $0$, as the Dirac measures near $q$ do in the paper.
\end{remark}

\section{Consistency of half the squared Fisher-Rao distance: proof of Theorem 4.6}\label{sec:S-theorem}

\begin{theorem}[Theorem 4.6]\label{thm:S-lcf}
For the Levi-Civita connection of the Fisher-Rao metric, the loss $\loss(q, x) = \tfrac{1}{2}d(q, x)^2$ is strictly $\Ff_0$-consistent for the Fr\'echet mean, which is defined, unique and interior for every $F \in \Ff_0$ and is the only solution of the barycenter equation $\E_F\log_q(X) = 0$.
\end{theorem}

Throughout, $\csimplex = \{p \in \R^n : p_i \ge 0,\ \sum_i p_i = 1\}$ is the closed simplex and $\Ff_0$ the class of Borel probability measures on the open simplex $\simplex$. The affinity $\rho(q, x) = \sum_i\sqrt{q_ix_i}$ is defined and continuous on $\csimplex \times \csimplex$, takes values in $[0, 1]$ by the Cauchy-Schwarz inequality, and equals $1$ exactly when $q = x$. The distance of Section~3 of the paper is $d(q, x) = 2\arccos\rho(q, x)$, so $\loss(q, x) = h(\rho(q, x))$ with $h(\rho) = 2(\arccos\rho)^2$. The function $h$ is continuously differentiable on $[0, 1]$, with $h'(\rho) = -4\arccos(\rho)/\sqrt{1 - \rho^2}$ for $\rho < 1$ and $h'(1) = -4$, since $(\arccos\rho)^2 = 2(1 - \rho) + O((1 - \rho)^2)$ as $\rho \to 1$. Writing $\rho = \cos\theta$ with $\theta \in [0, \pi/2]$, $h'(\rho) = -4\theta/\sin\theta$ decreases from $-4$ at $\theta = 0$ to $-2\pi$ at $\theta = \pi/2$, so
\begin{equation}\label{eq:S-hprime}
-2\pi \le h'(\rho) \le -4 \quad\text{for } \rho \in [0, 1].
\end{equation}
For $F \in \Ff_0$ let $R(q) = \E_F\,\loss(q, X) = \E_F\,h(\rho(q, X))$ be the expected loss, which is at most $\pi^2/2$. The Fr\'echet mean of $F$ is the minimizer of $\E_F\,d(q, X)^2 = 2R(q)$.

\begin{proof}
\emph{Step 1: $R$ on the closed simplex.} For each $x$, $q \mapsto h(\rho(q, x))$ is continuous on $\csimplex$ and bounded by $\pi^2/2$, so $R$ is continuous on $\csimplex$ by dominated convergence. Since $\csimplex$ is compact, $R$ attains its minimum on it.

\emph{Step 2: no minimizer on the boundary.} Let $q \in \csimplex \setminus \simplex$, let $I = \{i : q_i = 0\}$, which is nonempty, let $c = (1/n, \dots, 1/n)$ and $q(s) = (1 - s)q + sc$ for $s \in [0, 1]$, which lies in $\simplex$ for $s > 0$. For $x \in \simplex$ and $s \in (0, 1]$,
\[
\rho(q(s), x) = \sqrt{s/n}\,\sum_{i \in I}\sqrt{x_i} + \sum_{i \notin I}\sqrt{((1 - s)q_i + s/n)\,x_i},
\]
and its derivative with respect to $s$ is $A(s, x) + J(s, x)$ with
\[
A(s, x) = \frac{1}{2\sqrt{ns}}\sum_{i \in I}\sqrt{x_i} \ge 0, \qquad J(s, x) = \sum_{i \notin I}\frac{(1/n - q_i)\sqrt{x_i}}{2\sqrt{(1 - s)q_i + s/n}}.
\]
With $m = \min(\{q_i : i \notin I\} \cup \{1/n\}) > 0$ the denominators in $J$ are at least $2\sqrt m$, and $|1/n - q_i| \le 1$ and $\sqrt{x_i} \le 1$, so $|J(s, x)| \le K$ with $K = n/(2\sqrt m)$, a constant that depends on $q$ only. By \eqref{eq:S-hprime} and $A \ge 0$,
\[
\frac{d}{ds}\,h(\rho(q(s), x)) = h'(\rho(q(s), x))\,(A(s, x) + J(s, x)) \le -4A(s, x) + 2\pi K.
\]
For $0 < s_0 < s \le 1$ this derivative is bounded uniformly in $x$, so it can be integrated over $[s_0, s]$ under $\E_F$, and with $a = \E_F\sum_{i \in I}\sqrt{X_i}$,
\[
R(q(s)) - R(q(s_0)) \le -\frac{4a}{\sqrt n}\,(\sqrt s - \sqrt{s_0}) + 2\pi K(s - s_0).
\]
Here $a > 0$, since $X_i > 0$ almost surely under $F \in \Ff_0$. Letting $s_0 \to 0$ and using the continuity of $R$ from Step~1,
\[
R(q(s)) - R(q) \le -\frac{4a}{\sqrt n}\sqrt s + 2\pi Ks < 0 \quad\text{for } 0 < s < \left(\frac{2a}{\pi K\sqrt n}\right)^2.
\]
So $q$ is not a minimizer of $R$ on $\csimplex$, and every minimizer lies in $\simplex$.

\emph{Step 3: $R$ on the open simplex.} For $q, x \in \simplex$ the function $\loss(\cdot, x)$ is smooth near $q$, with gradient $\grad_q\loss(\cdot, x) = -\log_q(x)$ \citep[Proposition~3.10]{navasconvex}, a vector of length $d(q, x)$, and its Hessian at $q$ has the eigenvalue $1$ in the direction of $\log_q(x)$ and the eigenvalue $(\delta/2)\cot(\delta/2)$ on the orthogonal complement, $\delta = d(q, x)$, with the value $1$ at $\delta = 0$ \citep[Section~3.4]{navasconvex}; under the isometry $\Pi$ of Section~3 of the paper this is the Hessian of half the squared distance on a sphere of radius $2$. Since $(\delta/2)\cot(\delta/2) \le 1$,
\[
\Hess_q\loss(\cdot, x)(v, v) \ge \frac{d(q, x)}{2}\cot\frac{d(q, x)}{2}\;g_q(v, v), \qquad v \in T_q\simplex.
\]
Let $Q \subseteq \simplex$ be compact and $m_Q = \min\{q_i : q \in Q,\ i \le n\} > 0$. For $q \in Q$ and $x \in \csimplex$, $\rho(q, x) \ge \sqrt{m_Q}\sum_i\sqrt{x_i} \ge \sqrt{m_Q}$, so $d(q, x) \le \delta_Q$ with $\delta_Q = 2\arccos\sqrt{m_Q} < \pi$, and the eigenvalues of the Hessian lie in $[c_Q, 1]$ with $c_Q = (\delta_Q/2)\cot(\delta_Q/2) > 0$. The first and second derivatives of $\loss(\cdot, x)$ are thus bounded on $Q$ uniformly in $x$, so differentiation under $\E_F$ is allowed: $R$ is twice continuously differentiable on $\simplex$ with
\[
\grad_q R = -\E_F\log_q(X), \qquad \Hess_q R(v, v) \ge \E_F\left[\frac{d(q, X)}{2}\cot\frac{d(q, X)}{2}\right]g_q(v, v) \ge c_Q\,g_q(v, v)
\]
for $q \in Q$. At $q = \bary(F)$ the expectation in the middle is the constant $c_F$ of Section~7 of the paper. The critical points of $R$ are therefore exactly the solutions of the barycenter equation $\E_F\log_q(X) = 0$, and $\Hess_q R$ is positive definite at every $q \in \simplex$.

\emph{Step 4: at most one critical point.} The open simplex is geodesically convex: for $q, q' \in \simplex$ the points $\Pi(q)$ and $\Pi(q')$ lie in the open positive orthant of the sphere of radius $2$ at angular distance $d(q, q')/2 < \pi/2$, the shorter great-circle arc between them consists of positive combinations of the two points and so stays in the open orthant, and its preimage $\gamma: [0, 1] \to \simplex$ under $\Pi$ is the unique geodesic from $q$ to $q'$. Let $q$ be a critical point of $R$ and $q' \ne q$. The function $f(t) = R(\gamma(t))$ satisfies $f'(0) = g_q(\grad_q R, \dot\gamma(0)) = 0$ and $f''(t) = \Hess_{\gamma(t)}R(\dot\gamma(t), \dot\gamma(t)) > 0$, since $\gamma$ is a geodesic. So $f' > 0$ on $(0, 1]$ and $R(q') = f(1) > f(0) = R(q)$. A critical point is therefore the unique minimizer of $R$ on $\simplex$, and there is at most one.

\emph{Step 5: conclusion.} By Steps~1 and~2, $R$ has a minimizer $q^\ast$ in $\simplex$. It is a critical point, so $\E_F\log_{q^\ast}(X) = 0$. By Step~4 it is the only critical point, hence the only solution of the barycenter equation, and $R(q) > R(q^\ast)$ for every $q \in \simplex$ with $q \ne q^\ast$; by Step~2 the inequality also holds for $q \in \csimplex \setminus \simplex$. So the Fr\'echet mean $\bary(F) = q^\ast$ is defined, unique and interior for every $F \in \Ff_0$, and $\E_F\,\loss(\bary(F), X) < \E_F\,\loss(q, X)$ for $q \ne \bary(F)$, which is strict $\Ff_0$-consistency.
\end{proof}

\section{The cubic closedness defect at \texorpdfstring{$\alpha = \pm 1/2$}{alpha = 1/2 and -1/2}: numerical check of Proposition 5.6}\label{sec:S-numerics}

Proposition~5.6 of the paper states that for $n \ge 3$, $\alpha \notin \{-1, 0, 1\}$, $q \in \simplex$, $u, v, L \in T_q\simplex$ and $x = \exp_q(tL)$, the one-form $\omega_x = \musflat{\bar X^{(\alpha)}(\cdot, x)}$ satisfies
\begin{equation}\label{eq:S-cubic}
d\omega_x(u, v) = t^3\,\frac{(1 - \alpha^2)\alpha}{48}\,B_L(u, v) + O(t^4), \qquad B_L(u, v) = \sum_i r_i^2v_i\sum_j r_ju_j - \sum_i r_i^2u_i\sum_j r_jv_j,
\end{equation}
with $r_i = L_i/q_i$, and that $B_L \ne 0$ when three of the ratios $r_i$ are distinct. We compute $d\omega_x(u, v)$ for $n = 3$ directly from the geodesics, without the expansion, at $\alpha = 1/2$ and $\alpha = -1/2$, and at $\alpha = 0$, where the one-form is closed (Section~6 of the paper) and the computed value measures the numerical error.

The computation uses the mixture chart $p \mapsto (p_1, p_2)$, $p_3 = 1 - p_1 - p_2$, in which $g_{ij} = \delta_{ij}/p_i + 1/p_3$. The Christoffel symbols of $\bar\nabla^{(m)}$ vanish in this chart, and those of $\bar\nabla^{(e)}$, obtained by differentiating the exponential geodesics of Section~3 of the paper twice at $t = 0$, are $\Gamma^{(e)}_{ij,k} = g_{kl}\Gamma^{(e)l}_{ij} = -T_{ijk}$ with $T_{ijk} = \delta_{ijk}/p_i^2 - 1/p_3^2$ for $i, j, k \in \{1, 2\}$. By the definition of $\bar\nabla^{(\alpha)}$ as a convex combination, $\Gamma^{(\alpha)}_{ij,k} = -\frac{1 + \alpha}{2}\,T_{ijk}$. The geodesic equation $\ddot\gamma^k + \Gamma^{(\alpha)k}_{ij}\dot\gamma^i\dot\gamma^j = 0$ is integrated with an explicit Runge-Kutta method of order $8$ (DOP853 in SciPy) at relative tolerance $10^{-12}$ and absolute tolerance $10^{-14}$; this provides $\exp_q$. The residual $\bar X^{(\alpha)}(q, x) = \exp_q^{-1}(x)$ is obtained from projective flatness: by Equation~(36) of \citet{ayamari2015}, restated in Section~6 of the paper, the geodesic from $q$ to $x$ leaves $q$ in the direction of the projected cone residual $w = P_qX^{(\alpha)}(q, x)$, so $\bar X^{(\alpha)}(q, x) = sw$ for a scalar $s > 0$, and $s$ is the time at which the geodesic with initial velocity $w$ reaches $x$. It is found by bracketing and Brent's method on one coordinate, at tolerance $10^{-15}$, and the other coordinate is checked to agree with that of $x$ to $10^{-9}$. The one-form is evaluated as $\omega_x(w) = g_q(\bar X^{(\alpha)}(q, x), w)$. The coordinate vector fields $u$ and $v$ of the chart commute, so $d\omega_x(u, v) = \partial_u(\omega_x(v)) - \partial_v(\omega_x(u))$, and the two derivatives at $q$ are computed by central differences in $q$ with step $h$, once with $h = 4 \times 10^{-3}$ and once with $h = 2 \times 10^{-3}$, and combined by Richardson extrapolation, $\tfrac{4}{3}D_{h/2} - \tfrac{1}{3}D_h$, which removes the error term of order $h^2$. The script that produced the table and the figure is available on the author's webpage.

The base point and the vectors are
\[
q = (0.5, 0.3, 0.2), \qquad L = (0.3, -0.1, -0.2), \qquad u = (1, 0, -1), \qquad v = (0, 1, -1),
\]
so $u$ and $v$ are the coordinate vectors of the chart, and $r = (3/5, -1/3, -1)$ has three distinct ratios. The four sums in $B_L(u, v)$ are
\[
\begin{aligned}
\sum_i r_i^2v_i &= \tfrac{1}{9} - 1 = -\tfrac{8}{9}, &\qquad \sum_j r_ju_j &= \tfrac{3}{5} + 1 = \tfrac{8}{5},\\
\sum_i r_i^2u_i &= \tfrac{9}{25} - 1 = -\tfrac{16}{25}, &\qquad \sum_j r_jv_j &= -\tfrac{1}{3} + 1 = \tfrac{2}{3},
\end{aligned}
\]
so $B_L(u, v) = -\tfrac{64}{45} + \tfrac{32}{75} = -\tfrac{224}{225} = -0.9956$. At $\alpha = \pm 1/2$ the coefficient $(1 - \alpha^2)\alpha/48$ equals $\pm 1/128$, and the cubic term of \eqref{eq:S-cubic} is $\mp(224/225)\,t^3/128 = \mp 0.007778\,t^3$, the upper sign for $\alpha = 1/2$; at $t = 0.1$ it equals $\mp 7.778 \times 10^{-6}$.

\begin{table}[t]
\caption{Computed closedness defect $d\omega_x(u, v)$ at $q = (0.5, 0.3, 0.2)$ for $x = \exp_q(tL)$, and its ratio to the cubic term of \eqref{eq:S-cubic}. The last column is in units of $10^{-10}$.}\label{tab:closedness}
\begin{tabular}{lccccc}
\hline
& \multicolumn{2}{c}{$\alpha = 1/2$} & \multicolumn{2}{c}{$\alpha = -1/2$} & $\alpha = 0$ \\
$t$ & $d\omega_x(u, v)$ & ratio & $d\omega_x(u, v)$ & ratio & $d\omega_x(u, v) \times 10^{10}$ \\
\hline
$0.4$ & $-4.893 \times 10^{-4}$ & $0.983$ & $5.445 \times 10^{-4}$ & $1.094$ & $-6.06$ \\
$0.3$ & $-2.074 \times 10^{-4}$ & $0.988$ & $2.239 \times 10^{-4}$ & $1.066$ & $-5.37$ \\
$0.2$ & $-6.173 \times 10^{-5}$ & $0.992$ & $6.480 \times 10^{-5}$ & $1.041$ & $-4.48$ \\
$0.15$ & $-2.610 \times 10^{-5}$ & $0.994$ & $2.704 \times 10^{-5}$ & $1.030$ & $-4.00$ \\
$0.1$ & $-7.748 \times 10^{-6}$ & $0.996$ & $7.929 \times 10^{-6}$ & $1.019$ & $-3.49$ \\
$0.07$ & $-2.661 \times 10^{-6}$ & $0.997$ & $2.703 \times 10^{-6}$ & $1.013$ & $-3.16$ \\
$0.05$ & $-9.706 \times 10^{-7}$ & $0.998$ & $9.811 \times 10^{-7}$ & $1.009$ & $-2.94$ \\
$0.03$ & $-2.099 \times 10^{-7}$ & $1.000$ & $2.109 \times 10^{-7}$ & $1.004$ & $-2.70$ \\
\hline
\end{tabular}
\end{table}

Table~\ref{tab:closedness} lists $d\omega_x(u, v)$ for $t$ from $0.4$ to $0.03$, with $x = \exp_q(tL)$ computed by the geodesic solver, together with the ratio of the computed value to the cubic term. At $\alpha = 1/2$ the ratio increases from $0.983$ to $1.000$ as $t$ decreases, at $\alpha = -1/2$ it decreases from $1.094$ to $1.004$, and the computed defect changes sign with $\alpha$, as the coefficient $(1 - \alpha^2)\alpha$ does. At $\alpha = 0$ the computed values lie between $-6.1 \times 10^{-10}$ and $-2.7 \times 10^{-10}$.

\begin{figure}[t]
\centering
\includegraphics[width=0.78\textwidth]{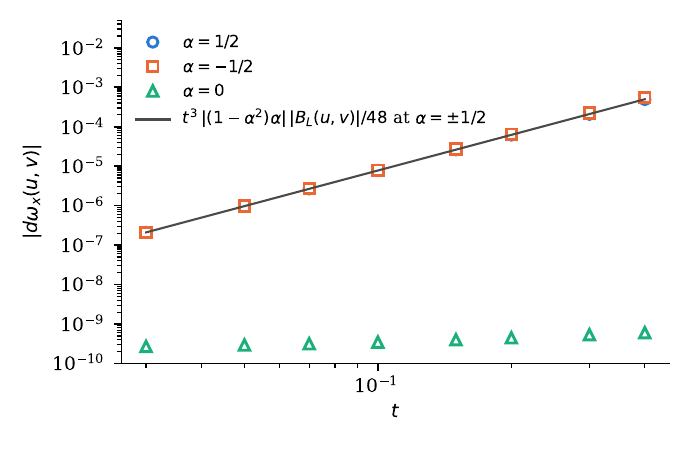}
\caption{Closedness defect $|d\omega_x(u, v)|$ of the one-form $\omega_x = \musflat{\bar X^{(\alpha)}(\cdot, x)}$ at $q = (0.5, 0.3, 0.2)$ against the geodesic parameter $t$, with $x = \exp_q(tL)$, on logarithmic axes. Markers: computed values at $\alpha = 1/2$ (circles), $\alpha = -1/2$ (squares) and $\alpha = 0$ (triangles). Line: the cubic term $t^3\,|(1 - \alpha^2)\alpha|\,|B_L(u, v)|/48$ of Proposition~5.6, which is the same at $\alpha = 1/2$ and $\alpha = -1/2$.}\label{fig:closedness}
\end{figure}

Figure~\ref{fig:closedness} plots the absolute values of the table on logarithmic axes. The markers are the computed defects $|d\omega_x(u, v)|$ at the eight values of $t$, circles for $\alpha = 1/2$, squares for $\alpha = -1/2$ and triangles for $\alpha = 0$. The grey line is the cubic term of Proposition~5.6, $t^3\,|(1 - \alpha^2)\alpha|\,|B_L(u, v)|/48 = (224/225)\,t^3/128$; the coefficient only changes sign between $\alpha = 1/2$ and $\alpha = -1/2$, so the line is the same for both. On logarithmic axes a cubic law is a straight line of slope $3$: the defect falls by three decades when $t$ falls by one.

Circles and squares lie on the line over the whole range. Between $t = 0.4$ and $t = 0.03$ the computed defect at $\alpha = 1/2$ decreases by the factor $2331$, against $(0.4/0.03)^3 = 2370$ for a cubic. This is the obstruction to closedness of the proposition, at both signs of $\alpha$. At the larger values of $t$ the markers deviate from the line, the circles below it and the squares above: by $1.7\%$ and $9.4\%$ at $t = 0.4$, by $0.4\%$ and $1.9\%$ at $t = 0.1$, and by less than $0.5\%$ at $t = 0.03$. The relative deviation decreases roughly in proportion to $t$, as the remainder $O(t^4)$ of the proposition requires. The sign of $d\omega_x(u, v)$, negative at $\alpha = 1/2$ and positive at $\alpha = -1/2$, is not visible in the figure, which plots absolute values; it is in Table~\ref{tab:closedness}. The triangles for $\alpha = 0$ stay between $3 \times 10^{-10}$ and $6 \times 10^{-10}$ for all $t$, nearly flat. This is the numerical floor of the procedure, the error of the central differences with step sizes $4 \times 10^{-3}$ and $2 \times 10^{-3}$ after Richardson extrapolation, and of the geodesic solver at tolerance $10^{-12}$. It lies more than two orders of magnitude below the smallest defect at $\alpha = \pm 1/2$, at $t = 0.03$, and six below the largest, at $t = 0.4$. At $\alpha = 0$ the one-form is closed, as Section~6 of the paper shows analytically, and the figure shows that the procedure separates a closed form from one with a cubic defect down to $t = 0.03$.

\end{document}